\documentclass[preprint,12pt]{elsart}

\usepackage{graphics}
\usepackage{amssymb}
\usepackage{lineno}

\usepackage{amsmath}
\usepackage{amsfonts}
\usepackage{amssymb}
\usepackage{amscd}
\usepackage{bbm}
\usepackage{bm}
\usepackage{placeins}

\usepackage{epstopdf,epsfig}
\newcommand{\RR}{\mathbbm{R}}

\def\v0{ {{\bf 0}} }

\def\vA{ {{\bf A}} }

\def\vf{ {{\bf f}} }

\def\vg{ {{\bf g}} }

\def\vM{ {{\bf M}} }

\def\vQ{ {{\bf Q}} }

\def\vu{ {{\bf u}} }

\def\vbeta{\mbox{\boldmath $\beta $}}

\def\vPsi{\mbox{\boldmath $\Psi $}}
\journal{Journal of Scientific Computing}

\begin{document}

\begin{frontmatter}



\title{On Spectral Approximations With Nonstandard Weight Functions and Their Implementations to Generalized Chaos Expansions}
%\tnoteref{t1}}
%\tnotetext[t1]{This research was supported by the ISRAEL SCIENCE
%FOUNDATION (grant No. 1364/04) and the UNITED STATES-ISRAEL
%BINATIONAL SCIENCE FOUNDATION (grant No. 2004099).}
%\author[tau]{S. Abarbanel} \ead{saul@post.tau.ac.il}
%\author[tau]{A.?Ditkowski\corref{cor}} \ead{adid@post.tau.ac.il}
%\cortext[cor]{Corresponding author}
%\address[tau]{School of Mathematical Sciences, Tel Aviv University, Tel Aviv 69978, Israel}
%
%\author{}
%
%\address{}

\author{A. Ditkowski}
\address{School of Mathematical Sciences, Tel Aviv University, Tel Aviv 69978, Israel}
\ead{adid@post.tau.ac.il}

\author{R. Kats}
\address{School of Mathematical Sciences, Tel Aviv University, Tel Aviv 69978, Israel}
\ead{ramikatz@mail.tau.ac.il}

\begin{abstract}
\indent   In this manuscript, we analyze the expansions of functions in orthogonal polynomials associated with a general weight function in a multidimensional setting. 
Such orthogonal polynomials can be obtained, e.g, by Gram-Schmidt orthogonalization. However, in most cases, they are not eigenfunctions of some singular Sturm-Liouville problem, as is the case for known polynomials, such as the Jacobi polynomials. Therefore, the standard convergence theorems do not apply.
Furthermore, since in general multidimensional cases the weight functions are not a tensor product of one-dimensional functions, the orthogonal polynomials are not a product of one-dimensional orthogonal polynomials, as well.

This work provides a way of estimating the convergence rate using a comparison lemma. We also present a spectrally convergent, multidimensional, integration method. Numerical examples demonstrate the efficacy of the proposed method. We also show that the use of non-standard weight functions can allow for efficient integration of singular functions. We demonstrate the use of this method to uncertainty quantification problem using Generalized Polynomial Chaos Expansions in the case of dependent random variables, as well.
\end{abstract}


\end{frontmatter}

%\linenumbers

\section{Introduction}\label{intro}
Orthogonal polynomials play a key role in approximation theory in weighted-$L^2$ spaces. Given a set $\Omega \subseteq \mathbb{R}^d$ and a weight function $w:\Omega \rightarrow [0,\infty)$, one can define a weighted-$L^2$ space, $L^2_w(\Omega)$, with an inner product induced by this weight function
\begin{equation}
\left<f,g\right>_w:=\int_{\Omega} f(x)g(x)w(x) dx
\end{equation}
In order to approximate a function $f \in L^2_w(\Omega)$, one can represent it as a series of polynomials that are orthogonal with respect to $\left<f,g\right>_w$. The required approximation is obtained by truncating the series.\bigskip \newline
In the case where the weight function is classical (s.a the constant or Gaussian weight functions) one can associate with the weight function a family of classical polynomials. For example, the orthogonal polynomials which are associated with the constant or Gaussian weight functions are the Legendre and Hermite polynomials, respectively. Moreover, these families of polynomials are eigenfunctions of known singular Sturm--Liouville problems. As a result, important properties such as the rate of convergence can be proved using the corresponding Sturm-Liouville operator. Moreover, it can be shown that the rate of convergence of such approximations depends on the smoothness of the approximated function. For more details see e.g. \cite{hesthaven2007spectral}, \cite{shen2011spectral}, \cite{funaro2008polynomial}, \cite{szeg1939orthogonal}, \cite{canuto1982approximation}. As a result of these properties, classical orthogonal polynomials have been proven to be a valuable tool in approximation theory and numerical analysis. An important application of such spectral approximations appears in the field of uncertainty quantification (UQ). This application motivated the current work and is described in the following discussion. However, we note that all the results described in this manuscript relate to questions in spectral approximations in general and are not limited to UQ.

Consider the PDE:
\begin{eqnarray}\label{1.10}
u_t(x,t,Z) &= &  {\cal L}(u)   \;\;,\;\;\;\;\;  x\in \Omega\;,
t \in[0,T]\;, Z\in \Omega_Z \subset \RR^d \nonumber \\
{\cal B}(u) &= & 0    \;,\;\;\;\;\;\;\;\;\;\;\;  x\in \partial \Omega\;,
t \in[0,T]\;, Z\in \Omega_Z \subset \RR^d \nonumber \\
u(t=0) &=& u_0 \;
\end{eqnarray}
where $Z = \left( z_1, \ldots,z_d\right)$ are random variables. A simple example of such a situation is the location of the transition layer in the solutions of Burgers equation, where the random variables are a perturbation in the boundary value. Another example is the heat equation in a domain whose properties are not exactly known. This problem can be modeled by a heat conduction coefficient which depends on random variables. In this case, the problem is to find the temperature or heat flux, at a given point, or the whole boundary. Equation \eqref{1.10} is a formulation of uncertainty quantification  problems which arise in many applications in science and engineering, See e.g. \cite{xiu2010numerical}, \cite{xiu2003new}, \cite{xiu2005high}, \cite{xiu2003modeling}, \cite{xiu2009fast},   \cite{xiu2006numerical}, \cite{gottlieb2008galerkin}, \cite{xiu2009fast}, \cite{lin2007random}.

Several methods were developed to analyze these problems, among which are Monte Carlo and Sampling-Based Methods, Perturbation Methods, Moment Equations and Generalized Polynomial Chaos.  In the GPC method the solution to \eqref{1.10}, $u(x,t,Z) $, which is an unknown random variable is expressed as a function of the random variables $Z$, whose distributions are known. Thus, $u(x,t,Z) $ can be expanded in orthogonal polynomials in $Z$. If $u(x,t,Z)$ is smoothly dependent on $Z$, and the number of $z_j$ is small, typically less than six, the GPC expansion is highly efficient due to a spectral rate of convergence, see e.g. \cite{xiu2010numerical}.

The standard assumption is that the random variables are  mutually independent, see e.g. \cite{xiu2010numerical}. More generally, the classical GPC theory in the multi-dimensional setting relies on the fact that the multivariate weight function can be decomposed into a product of one-dimensional weight functions
\begin{equation}\nonumber
w(Z) = \prod_{j=1}^d w_j(z_j)
\end{equation}
As a result, the integral of an integrable function $f:\Omega_Z \rightarrow \mathbb{R}$ can be calculated by applying Fubini's theorem
\begin{equation}\nonumber
\int_{\Omega_Z}  f(Z) w(Z)dZ =  \int w_1(z_1)\int w_2(z_2) \dots \int w_d(z_d)f(Z)dz_d \dots dz_1  \;.
\end{equation}
For each inner product
\begin{equation}\nonumber
\left(u,v \right )_{w_j} =  \int u v w_j(z_j ) dZ_j  
\end{equation}
a set of one-dimensional orthogonal polynomials can be associated and the basis for the expansion in the whole $\Omega_Z$ is a tensor product of these one-dimensional polynomials. Important properties, such as the rate of convergence and cubature integration formulas, are governed by the properties of the 1-D expansions.

This assumption makes the  expansion and its analysis much simpler as the  polynomial basis is a tensor product of one-dimensional  orthogonal polynomials. There are situations, however, in which the random variables are not independent. In such cases, the decomposition of the multivariate weight function can not be carried out. An example of such a situation is temperature and relative humidity. As the temperature changes the relative humidity changes as well. Drying processes  are strongly dependent on both the temperature and the relative humidity, see \cite{mujumdar2014handbook}.

There are two common approaches to address the dependence problem. The first one is to transform the problem into an independent variables problem using, for example, the  Rosenblatt   transformation, \cite{rosenblatt1952remarks}.  This transformation, however, may be ill-conditioned. The second approach is to 
consider the problem in an epistemic uncertainty framework, see e.g.  \cite{jakeman2010numerical}, \cite{chen2013flexible}, \cite{li2014upper}.
Since the exact distributions are not used in this approach, the convergence rate may be slow.

In the case of dependent variables, $w(Z)$ is indeed a genuine $d$ dimensional function which cannot be decomposed into a product of one-dimensional functions. Though a set of orthogonal polynomials is  defined on $\Omega_Z$ and they can be found, e.g. using the Gram--Schmidt process, there is no general theory of expansions in orthogonal polynomials in the multidimensional case, see e.g. \cite{xu2004lecture}.

Even in the 1-D case, it can be difficult to evaluate the rate of convergence for a general weight function,  $W(Z)$, since there is no singular Sturm-Liouville problem associated with it.

\textbf{Remark}: As we've mentioned before, these are questions in spectral approximations and are not limited to GPC expansions.

In this paper, we analyze the expansion of functions in the set of orthogonal polynomials associated with $w(Z)$ in $\Omega_Z$. We evaluate the convergence rate using a comparison to a known distribution. We also propose a  heuristic way  to derive cubature integration formulas. The paper is constructed as follows: in Section \ref{section:convergence_rate} we analyze the convergence rate of the orthogonal polynomial expansion. We present a Lemma which enables to determine the convergence rate by comparing it to an expansion with a known convergence rate. In Section \ref{section:numerical integration} we present a spectrally accurate integration method in multidimensional domains. In Section \ref{Example_of_implementation_to_GPC} we present the implementation of the proposed method in a GPC setting. Concluding remarks are presented in Section \ref{section:conclusions}.

\section{Convergence Rate}\label{section:convergence_rate}

Let
\begin{equation}\label{2.10}
f(Z) = \sum_{|\nu|=0}^\infty \hat{f}(\nu) P_{\nu}(Z) 
\end{equation}
where $\nu = (\nu_1,\ldots\nu_d)$, $\nu_j \ge 0$ is a multi--index and $\left\{P_{\nu}(Z)\right\}_{|\nu|=0}^{\infty}$ is a set of orthonormal polynomials.

The question is how fast $\hat{f}(\nu) \longrightarrow 0$ as $|\nu|\longrightarrow \infty$?. For the 1-D classical polynomials, s.a. Jacobi polynomials, there is an underlying singular  Sturm--Liouville problem, and the spectral convergence results from it. In the multidimensional case or in the case of non-standard weight functions, there is no similar theory, see e.g. \cite{xu2004lecture}. Therefore, we analyze the rate of convergence by comparison to known weight functions.

\bigskip

\subsection{Main Lemma}\label{main_theorem}
\textbf{Lemma:} Let $V$ be an infinite dimension linear vector space  and ${\cal{H}}_1$ and ${\cal{H}}_2$ be two separable Hilbert spaces defined on  $V$ with inner products $\left < \cdot,\cdot\right >_1$ and $\left < \cdot,\cdot\right >_2$, respectively.  Let $\{\phi_j^{k}\}_{j \geq 0}$ be an orthogonal basis of ${\cal{H}}_k$ with respect to $\left < \cdot,\cdot\right >_k$ and suppose that  
\begin{equation}\label{main_lemma_assumption_1} 
 \ \  {\rm span}\{\phi_j^{1}\}_{0\leq j \leq N}={\rm span}\{\phi_j^{2}\}_{0\leq j \leq N}=V_N ,\enskip \forall  N \in \mathbb{N}
\end{equation}
and
\begin{equation}\label{main_lemma_assumption_2} 
 \ \ \left\|  x\right\|_2 \, \leq \, C \left\|  x\right\|_1 \enskip \forall x \in V
\end{equation}
where $C>0 $ is a constant that is independent of $x$ .

Let $P_N^{k}$ be the projection operators onto $ {\rm span}\{\phi_j^{k}\}_{0\leq j \leq N}$ with respect to the inner product $\left < \cdot,\cdot\right >_k$, $k=1,2$ and $Q_N^{k}= I - P_N^{k}$ their complementary projection  operators.

Then, For all $ x \in V$,
\begin{equation}\label{main_lemma_1} 
\left\| Q_N^{2}x  \right\|_2\, \leq \, C \left\| Q_N^{1}x\right\|_1 
\end{equation}
and
\begin{equation}\label{main_lemma_2} 
\left\| P_N^{2}x  - P_N^{1}x\right\|_2\, \leq \, C \left\| Q_N^{1}x\right\|_1 
\end{equation}

\textbf{Proof:} Note that $x=P_N^{1}x+Q_N^{1}x$ in the sense that  $||x-\left( P_N^{1}x+Q_N^{1}x \right ) ||_1 = 0$. Similarly, $x=P_N^{2}x+Q_N^{2}x$ in the sense that  $||x-\left( P_N^{2}x+Q_N^{2}x \right ) ||_2 = 0$. Due to \eqref{main_lemma_assumption_2} $||x-\left( P_N^{1}x+Q_N^{1}x \right ) ||_2 = 0$ as well.

Let $x \in V$. Since $x=P_N^{1}x+Q_N^{1}x=P_N^{2}x+Q_N^{2}x$ in $||\cdot||_2$ sense, by using \eqref{main_lemma_assumption_2} we obtain the following estimate:
\begin{equation}
C^2 \left\|Q_N^{1}x\right\|_1^{2}\, \ge\, \left\|Q_N^{1}x\right\|_2^{2}\, = \, \left\|\left(P_N^{2}x-P_N^{1}x\right)+Q_N^{2}x\right\|_2^{2} 
\end{equation}
Since $P_N^{2}x-P_N^{1}x\in V_N$ and $Q_N^{2}x \in V_N^{\perp,2}$ 
\begin{equation}
C^2 \left\|Q_N^{1}x\right\|_1^{2}\, \ge\, \left\| P_N^{2}x-P_N^{1}x\right\|_2^{2} +\left\|Q_N^{2}x\right\|_2^{2} 
\end{equation}
which leads to \eqref{main_lemma_1} and \eqref{main_lemma_2}.

{\bf Remark}: In the case of comparing two sets of orthogonal polynomials, condition \eqref{main_lemma_assumption_1} is trivially satisfied. In this case, the conclusion  from the lemma is that the convergence rate of the expansion in polynomial derived from bounded  $W(Z)$  is, at least,  of the same order as the expansion in Legendre  polynomials.  

\subsection{Numerical examples}

In this section, we demonstrate the claims of the Lemma using several numerical examples.

\begin{figure}[h]
\begin{center}
\begin{tabular}{lllll}
a:Legendre&b:Gram Schmidt\\
\includegraphics[width=0.5\textwidth]{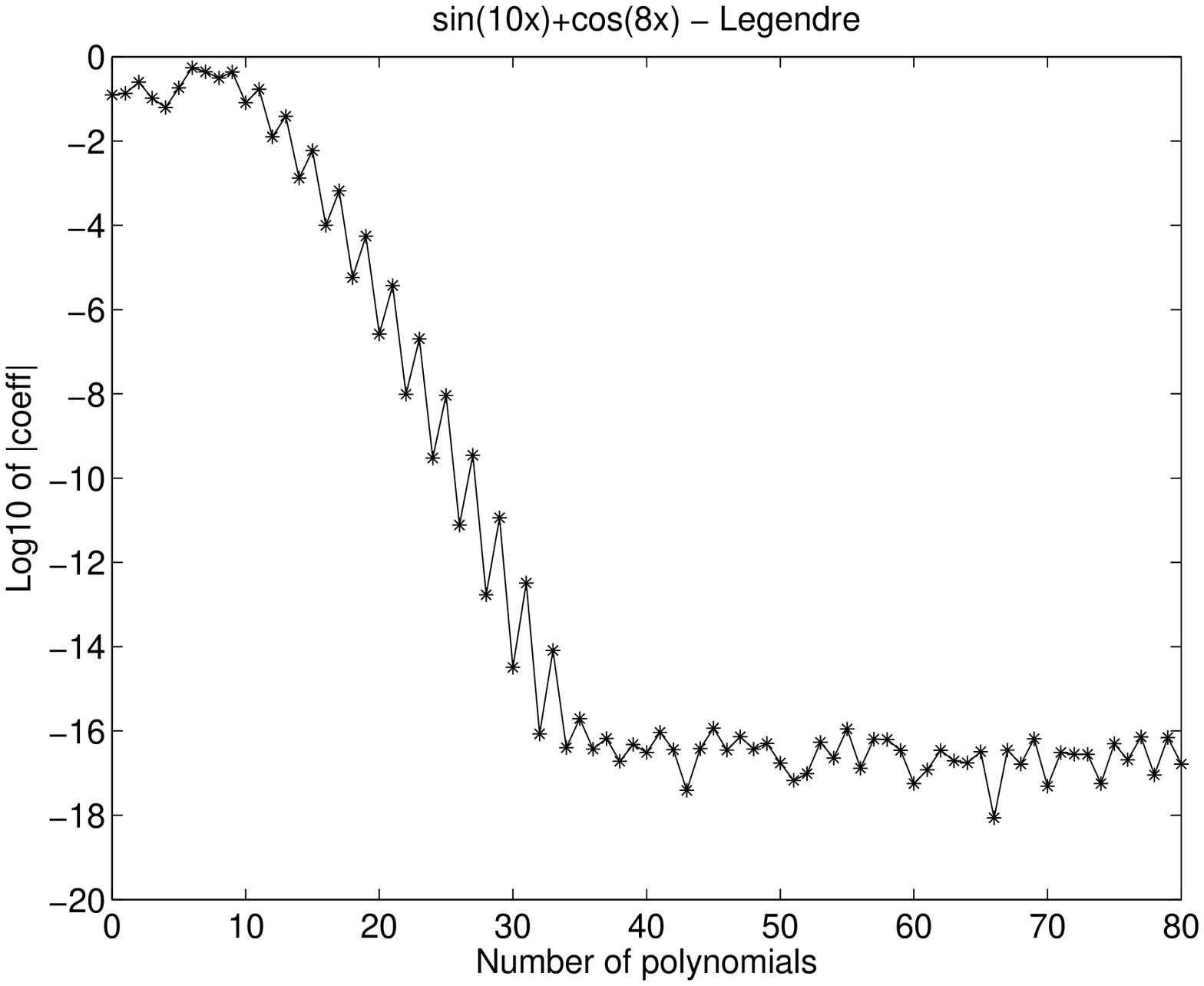}&
\includegraphics[width=0.5\textwidth]{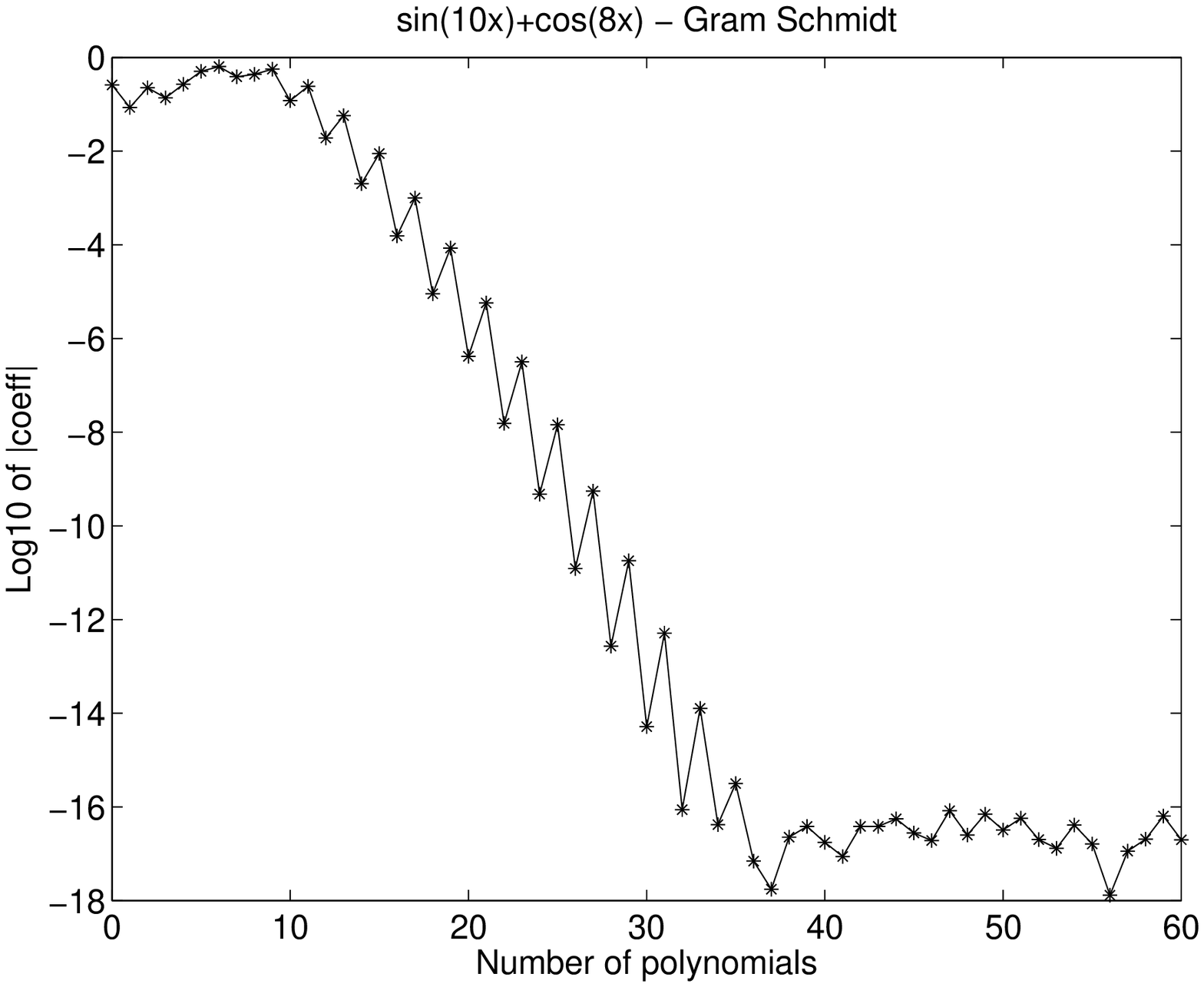}
\end{tabular}
\end{center}
\caption{Comparison between the coefficients of the expansions of $f(x)=\sin(10x)+\cos(8x)$ in a: Legendre polynomials and b: the orthogonal polynomials associated with $w_2(x)=\frac{3}{4}\left(1-x^2\right)$.
} \label{fig:fig_num_exam_1.1}
\end{figure}

\begin{figure}[h]
\begin{center}
\begin{tabular}{lllll}
a:Legendre&b:Gram Schmidt\\
\includegraphics[width=0.5\textwidth]{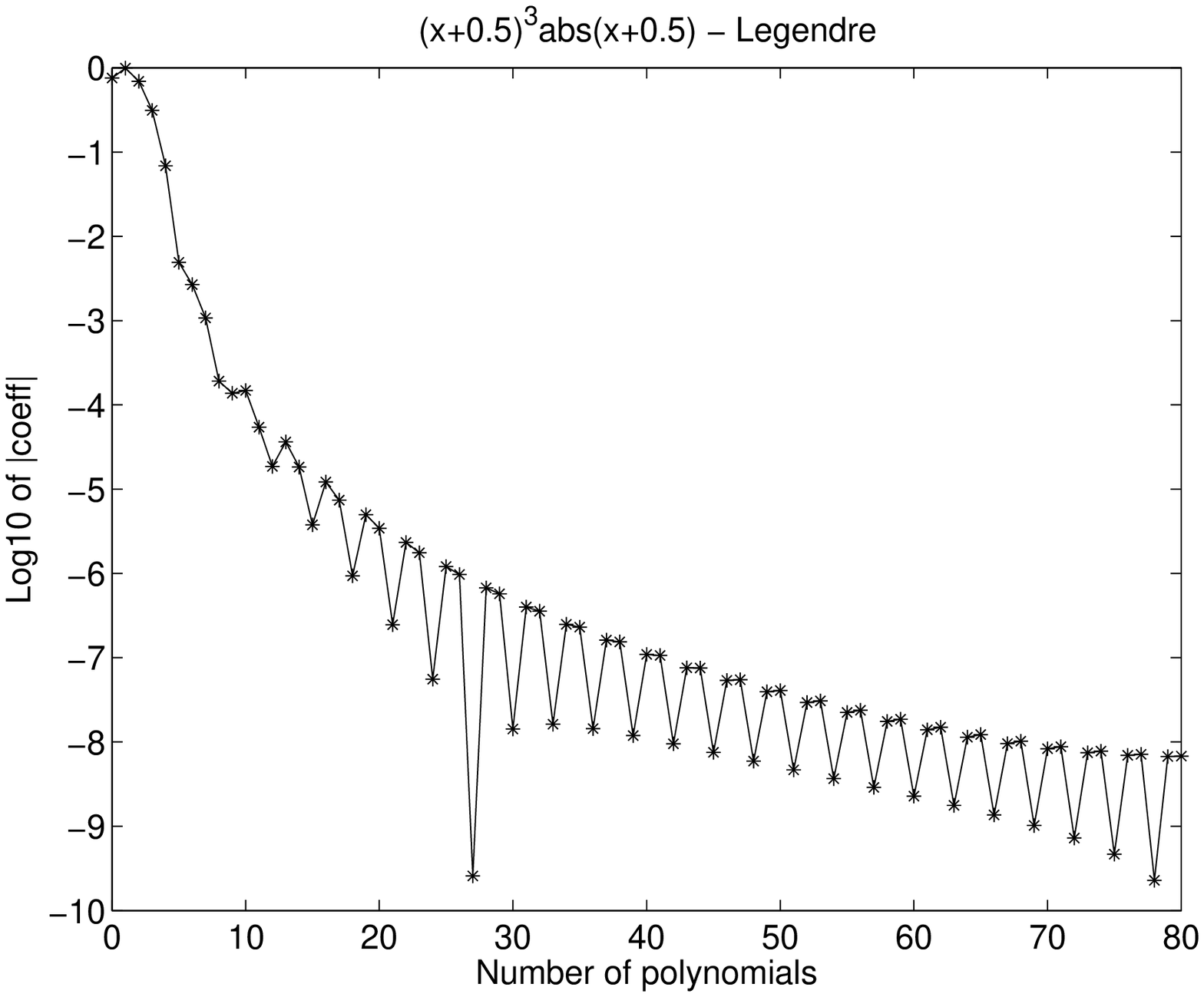}&
\includegraphics[width=0.5\textwidth]{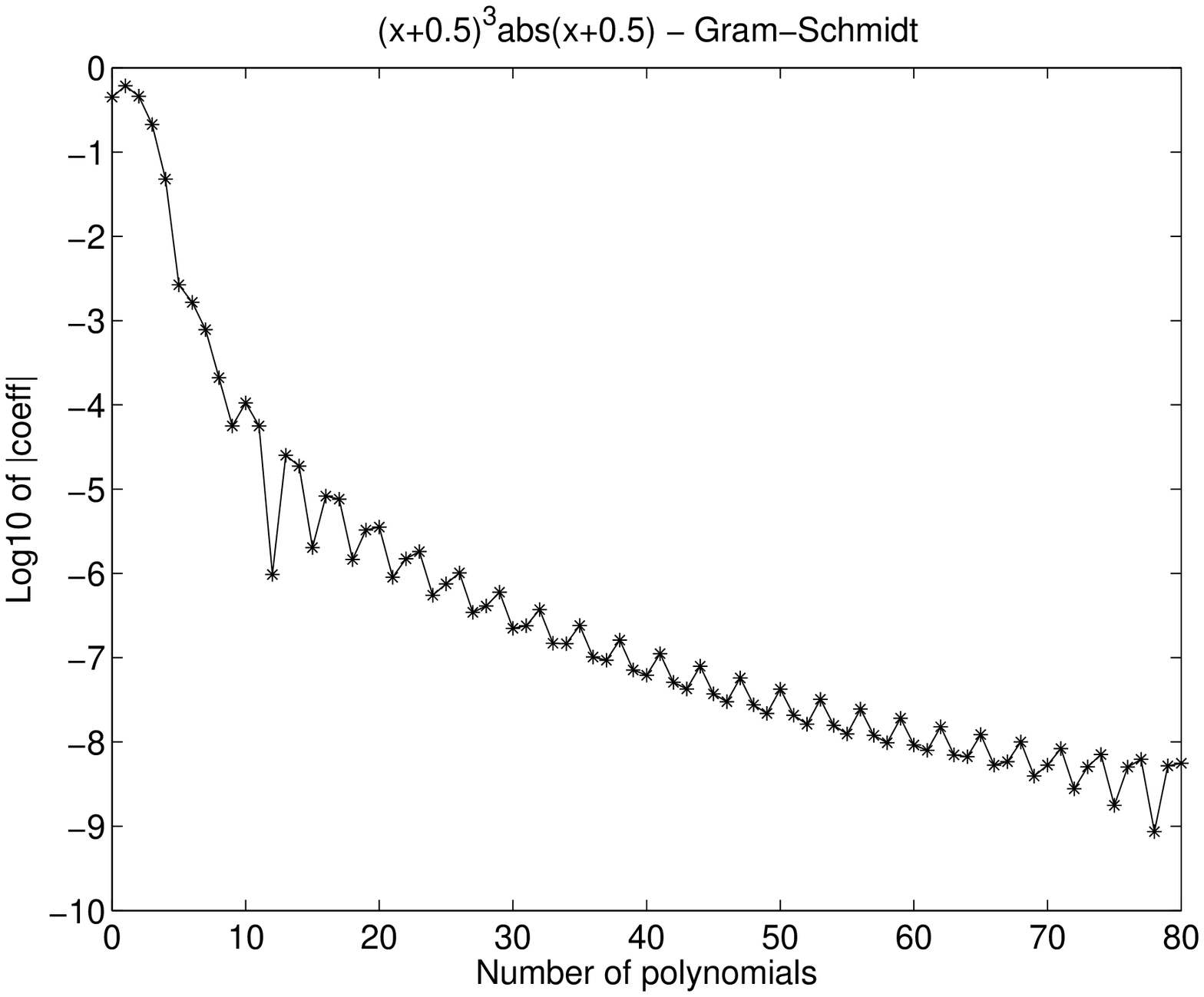}
\end{tabular}
\end{center}
\caption{Comparison between the coefficients of the expansions of $ g(x)=(x+ \frac{1}{2})^3 \left | x+ \frac{1}{2} \right |$ in a: Legendre polynomials and b: the orthogonal polynomials associated with  $w_2(x)=\frac{3}{4}\left(1-x^2\right)$.
} \label{fig:fig_num_exam_1.2}
\end{figure}

\begin{figure}[h]
\begin{center}
\begin{tabular}{lllll}
a:Legendre&b:Gram Schmidt\\
\includegraphics[width=0.5\textwidth]{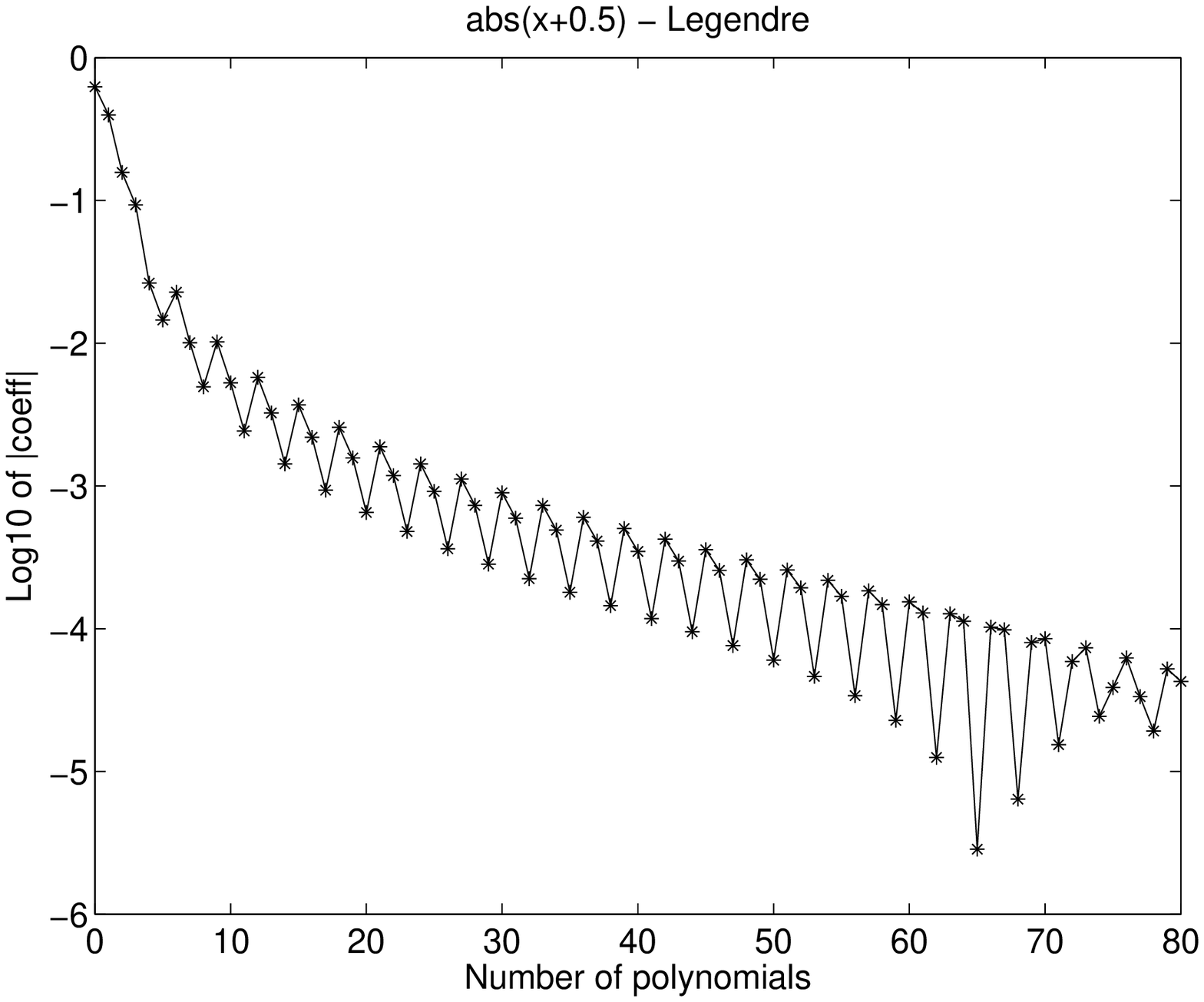}&
\includegraphics[width=0.5\textwidth]{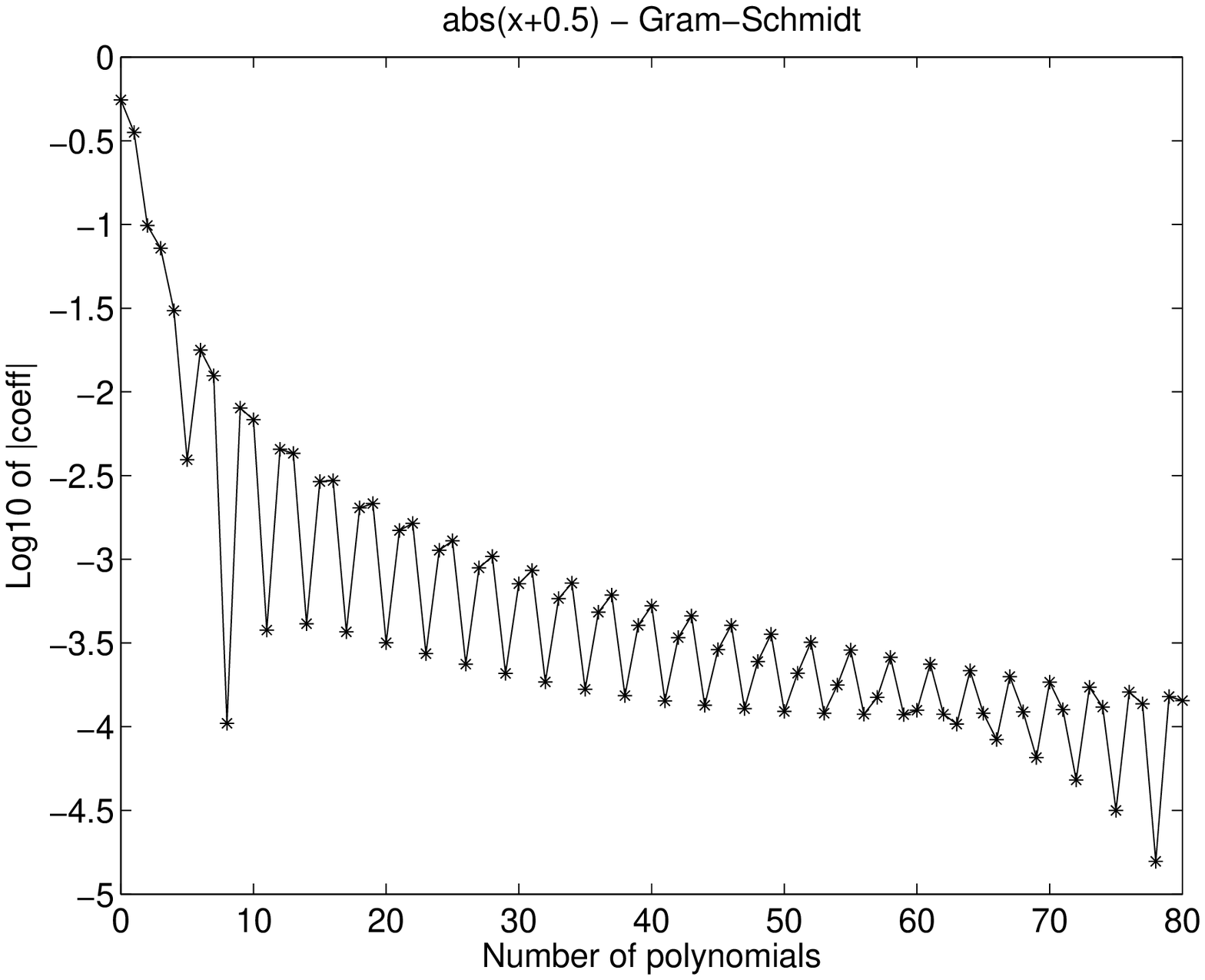}
\end{tabular}
\end{center}
\caption{Comparison between the coefficients of the expansions of $ h(x)=\left | x+ \frac{1}{2} \right |$ in a: Legendre polynomials and b: the orthogonal polynomials associated with  $w_2(x)=\frac{3}{4}\left(1-x^2\right)$.
} \label{fig:fig_num_exam_1.3}
\end{figure}

{\bf Example 1}: 
In this example we compare the rate of decay of the projection coefficients of the functions $ f(x)=\sin(10x)+\cos(8x)$, $ g(x)=(x+ {1}/{2})^3 \left | x+ {1}/{2} \right |$ and $ h(x)=\left | x+ {1}/{2} \right |$. \newline \\
We compare the expansions of this function with respect to Legendre-type polynomials ($w_1(x)={1}/{2}$) and non-standard orthogonal polynomials ($w_2(x)={3}/{4}\left(1-x^2\right)$). The latter were obtained via Gram-Schmidt orthogonalization. The plots of the logarithm of the coefficients are presented in Figures \ref{fig:fig_num_exam_1.1}, \ref{fig:fig_num_exam_1.2}, \ref{fig:fig_num_exam_1.3} ,respectively. Since $f(x)$ is analytic, $g(x) \in C^3[-1,1]$ and $h(x) \in C[-1,1]$ the decay rates of their coefficients are significantly different. However the comparison between the Legendre-type expansion and the expansion in the polynomials associated with $w_2(x)$ indicates that the decay rates are very much alike.   

We also note that the norm that is induced by $w_2(x)$ is dominated by the norm induced by $w_1(x)$ and, therefore, the conditions of the lemma are satisfied.

%$Cos(8x) + Sin(10y)$; Legendre$ (w=1/2)$, and $w=\frac{3}{4}(1-x^2)$

\begin{figure}[h]
\begin{center}
\begin{tabular}{lllll}
a:Chebychev&b:Gram Schmidt\\
\includegraphics[width=0.5\textwidth]{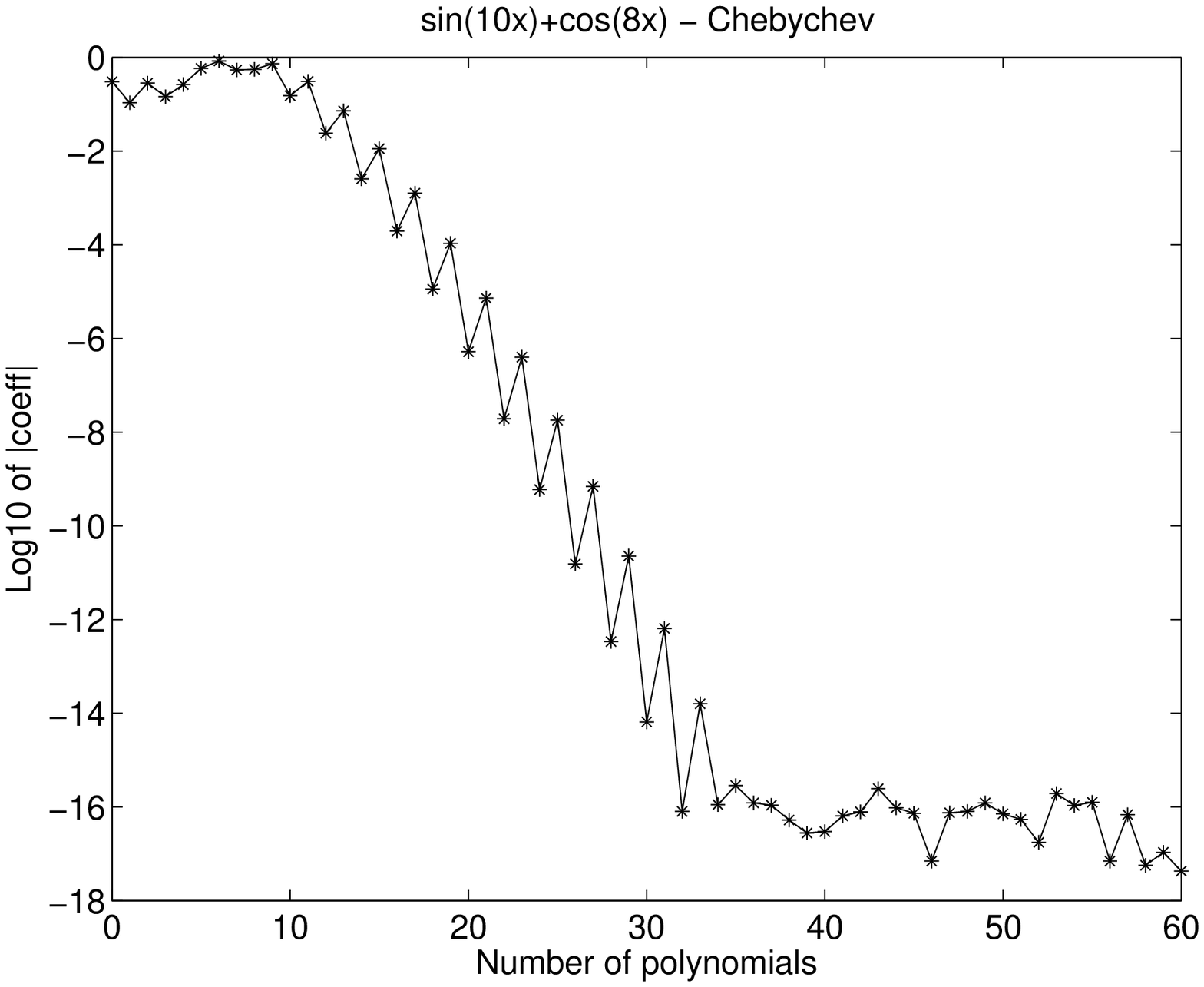}&
\includegraphics[width=0.5\textwidth]{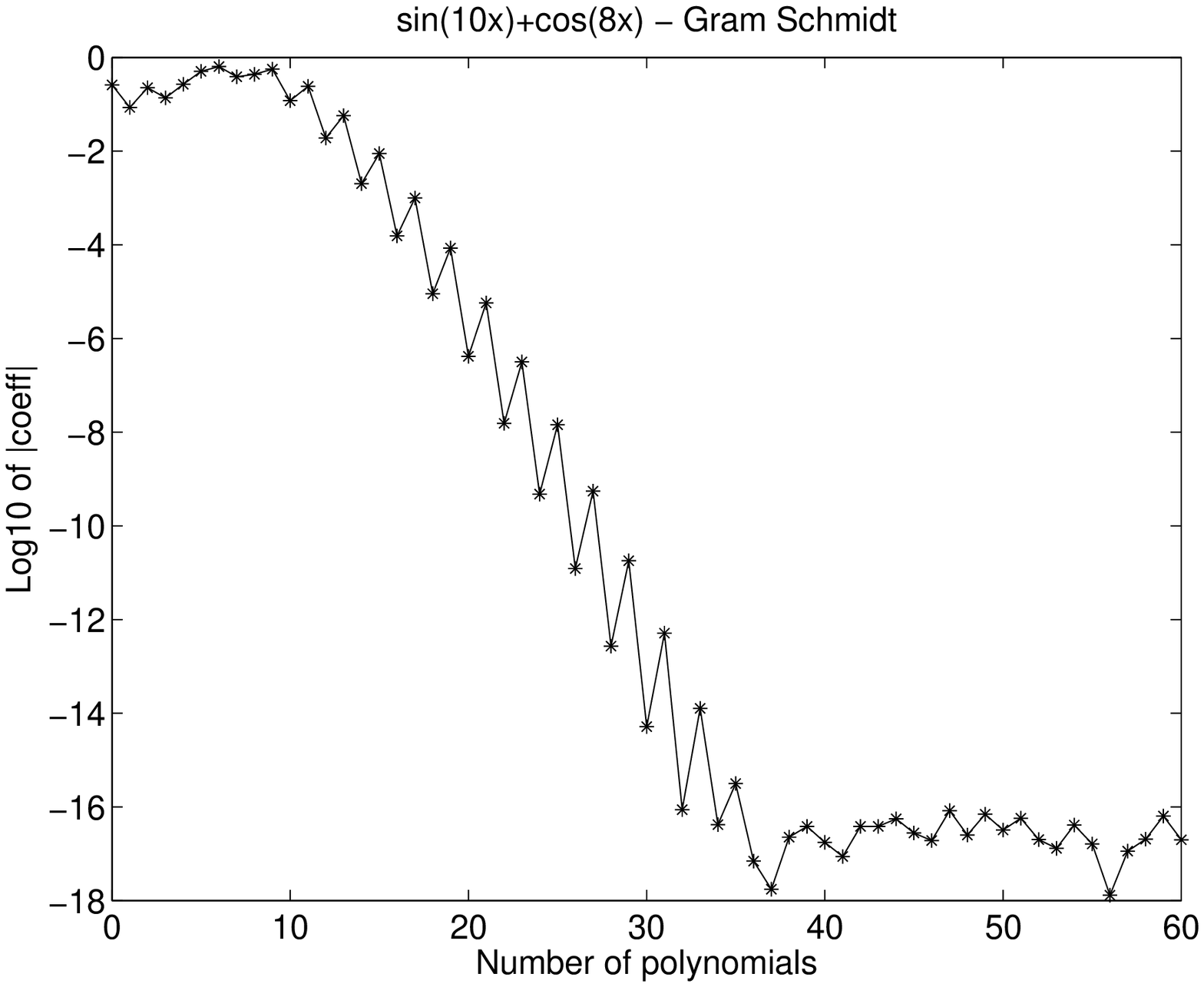}
\end{tabular}
\end{center}
\caption{Comparison between the coefficients of the expansions of $f(x)=\sin(10x)+\cos(8x)$ in a: Chebychev polynomials and b: the orthogonal polynomials associated with  $w_2(x)={1}/{(2\sqrt{1-|x|})}$.
} \label{fig:fig_num_exam_2.1}
\end{figure}

\begin{figure}[h]
\begin{center}
\begin{tabular}{lllll}
a:Chebychev&b:Gram Schmidt\\
\includegraphics[width=0.5\textwidth]{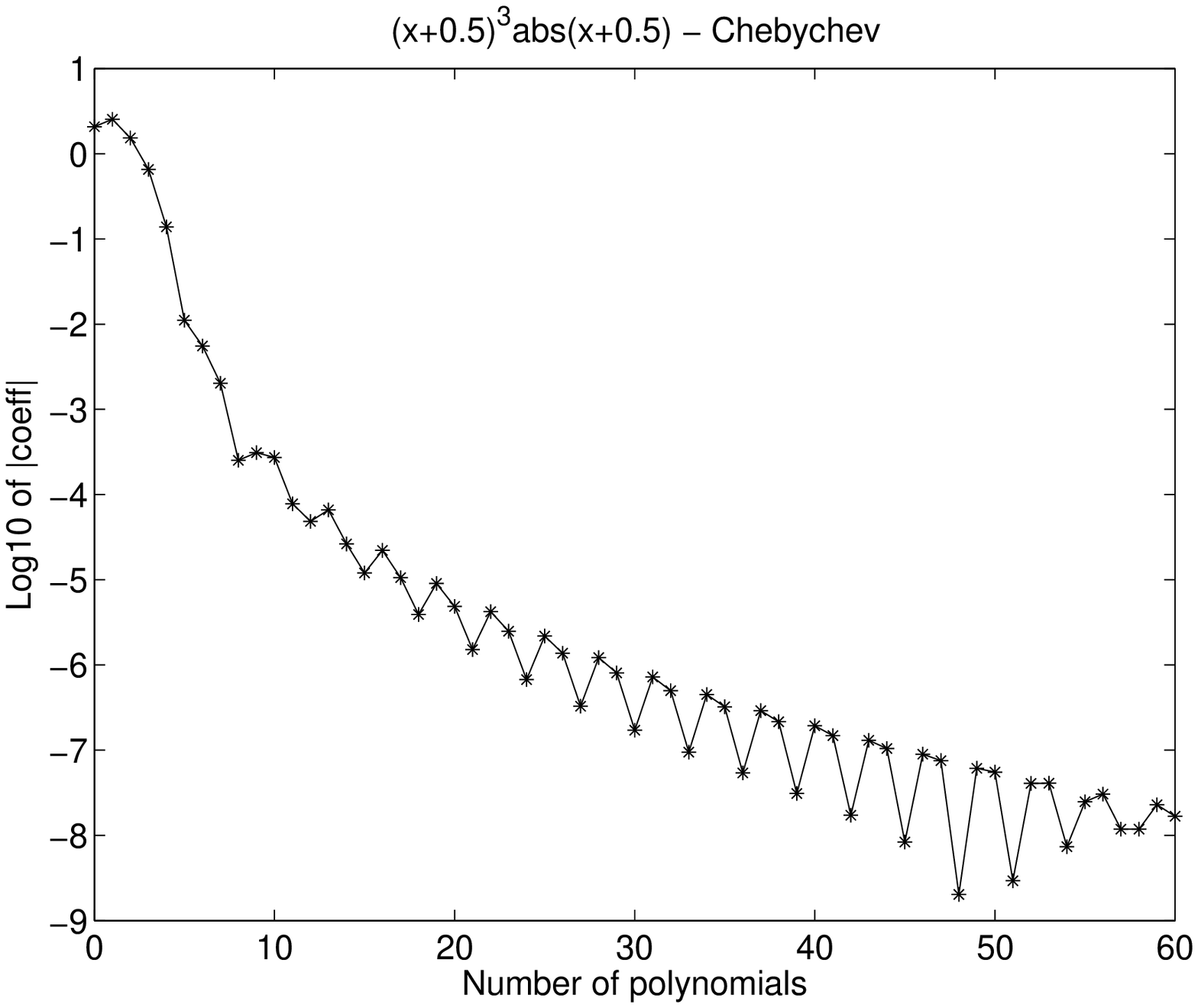}&
\includegraphics[width=0.5\textwidth]{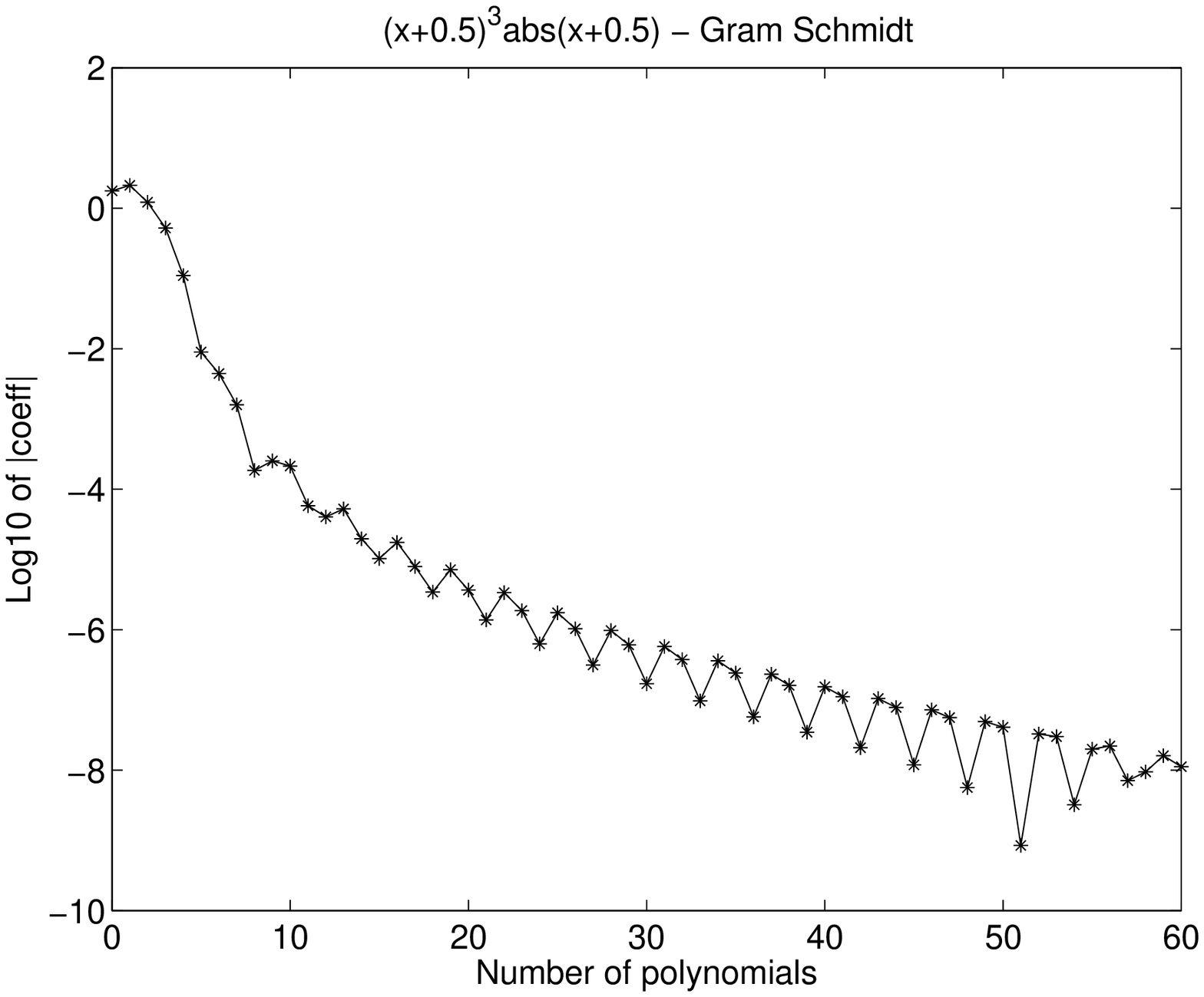}
\end{tabular}
\end{center}
\caption{Comparison between the coefficients of the expansions of $ g(x)=(x+ \frac{1}{2})^3 \left | x+ \frac{1}{2} \right |$ in a: Chebychev polynomials and b: the orthogonal polynomials associated with  $w_2(x)={1}/{(2\sqrt{1-|x|})}$.
} \label{fig:fig_num_exam_2.2}
\end{figure}

\begin{figure}[h]
\begin{center}
\begin{tabular}{lllll}
a:Chebychev&b:Gram Schmidt\\
\includegraphics[width=0.5\textwidth]{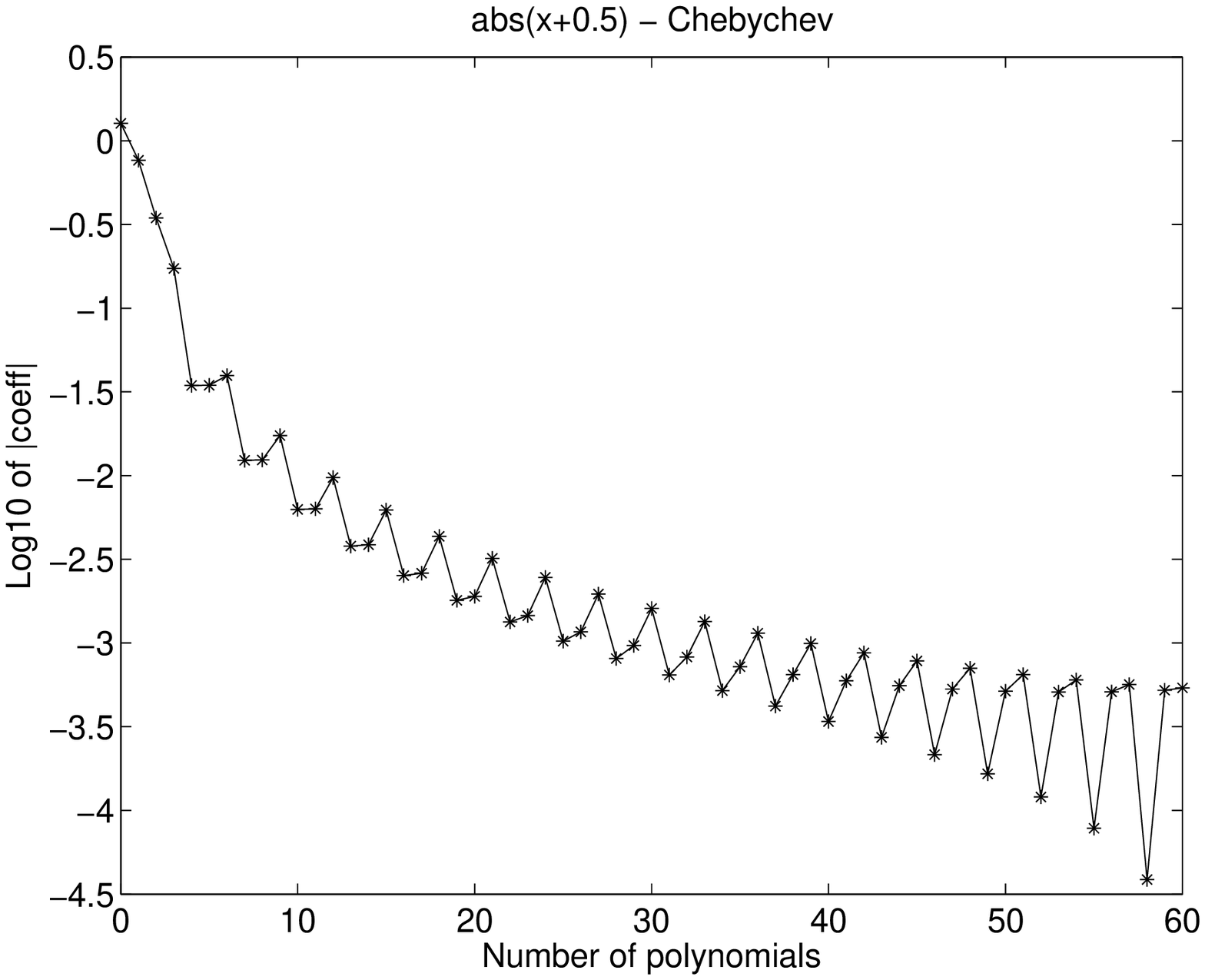}&
\includegraphics[width=0.5\textwidth]{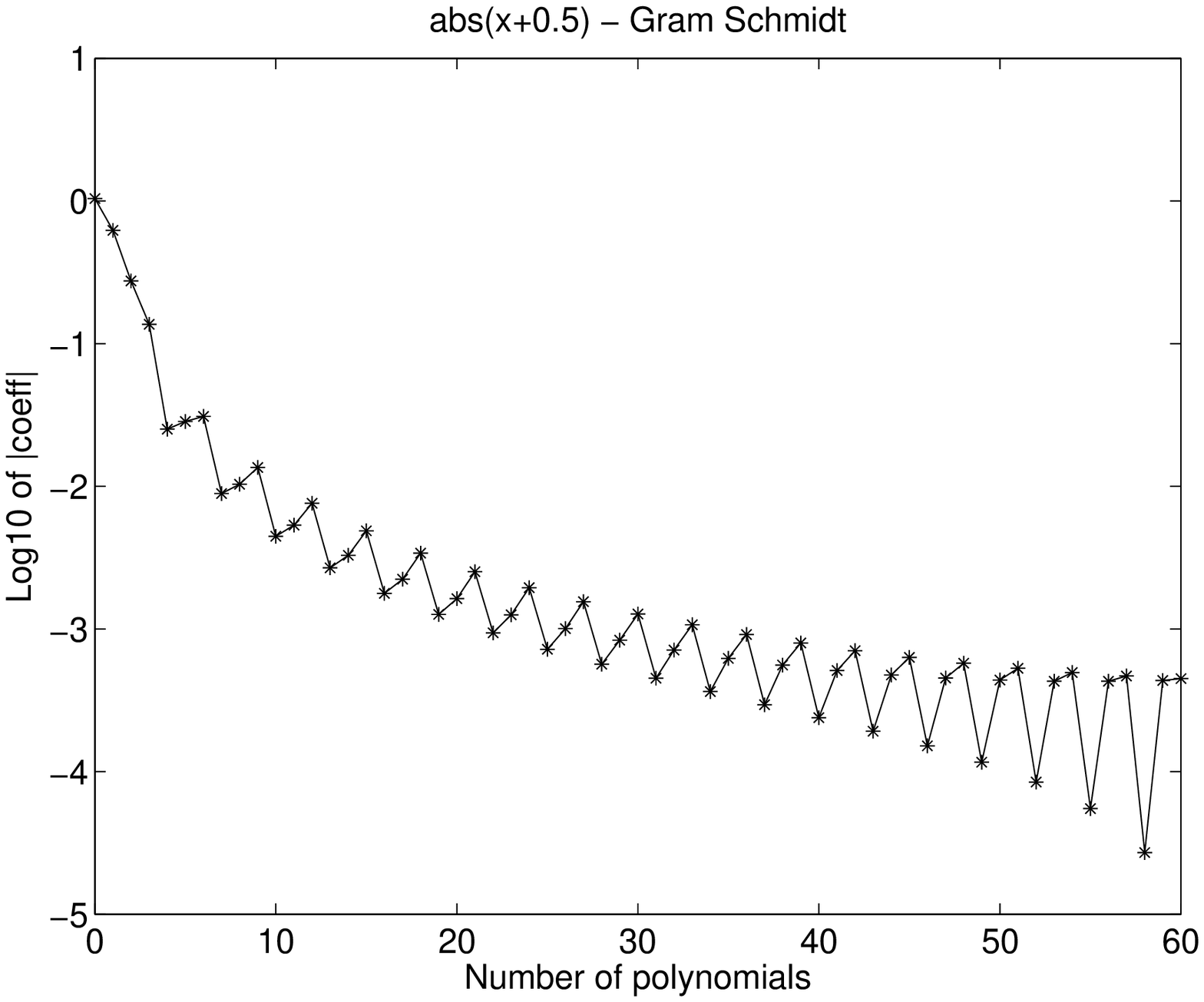}
\end{tabular}
\end{center}
\caption{Comparison between the coefficients of the expansions of $ h(x)=\left | x+ \frac{1}{2} \right |$ in a: Chebychev polynomials and b: the orthogonal polynomials associated with  $w_2(x)={1}/{(2\sqrt{1-|x|})}$.
} \label{fig:fig_num_exam_2.3}
\end{figure}

{\bf Example 2}: In this example we compare the rate of decay of the projection coefficients of the same functions,  $ f(x)=\sin(10x)+\cos(8x)$, $ g(x)=(x+ {1}/{2})^3 \left | x+ {1}/{2} \right |$ and $ h(x)=\left | x+ {1}/{2} \right |$. \newline \\ 
We compare the expansions of this function with respect to Chebychev polynomials ($w_1(x)={1}/{\sqrt{1-x^2}}$) and non-standard orthogonal polynomials ($w_2(x)={1}/{(2\sqrt{1-|x|})}$) which were obtained via Gram-Schmidt orthogonalization. The plots of the logarithm of the coefficients are presented in Figures \ref{fig:fig_num_exam_2.1}, \ref{fig:fig_num_exam_2.2}, \ref{fig:fig_num_exam_2.3}, respectively. As in the previous example, the decay rates are almost identical. Here as well, the norm that is induced by $w_2(x)$ is dominated by the norm induced by $w_1(x)$. 

%$Cos(8x) + Sin(10y)$; Legendre$ (w=1/2)$, and $w=\frac{3}{4}(1-x^2)$

\begin{figure}[h]
\begin{center}
\begin{tabular}{lllll}
a:Legendre&b:Gram Schmidt\\
\includegraphics[width=0.5\textwidth]{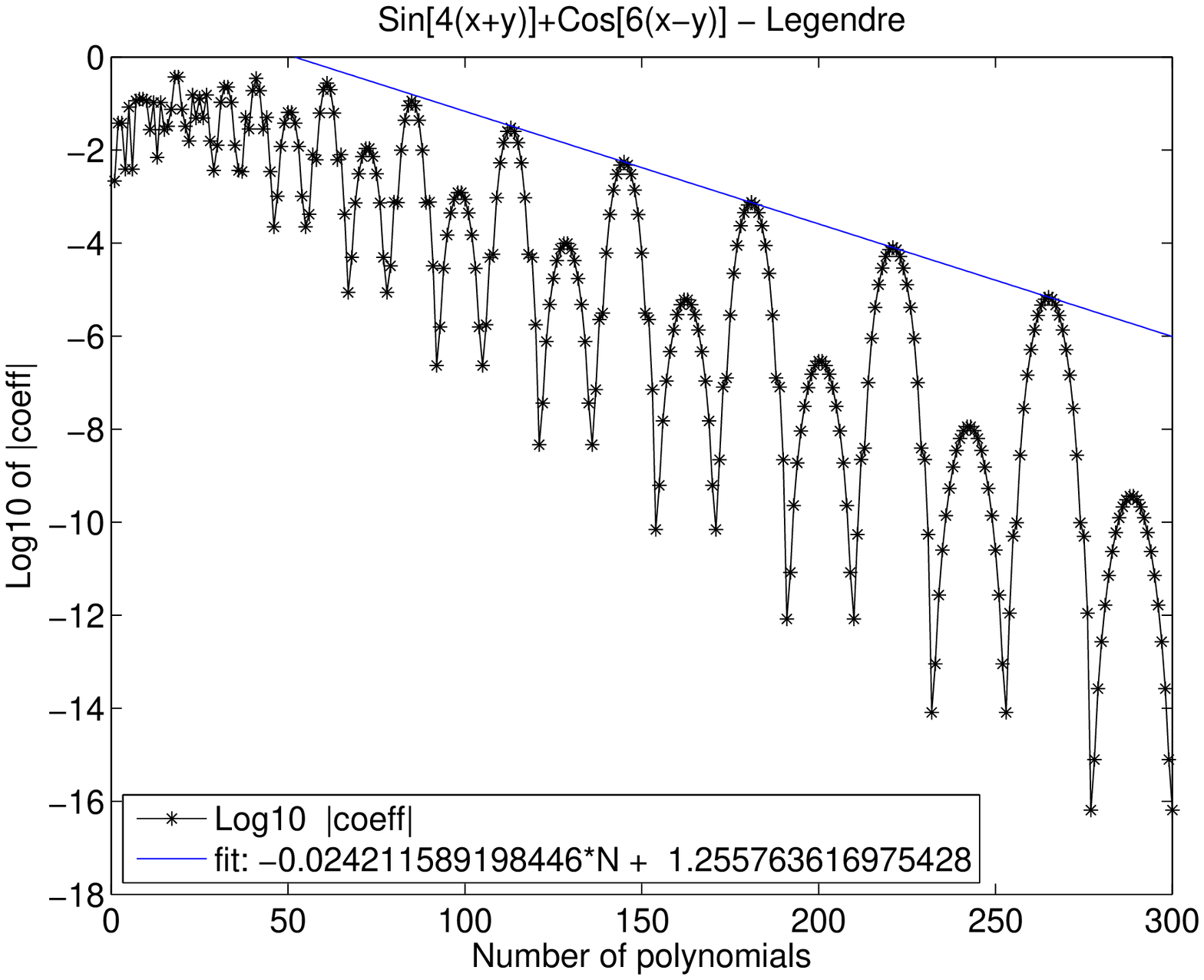}&
\includegraphics[width=0.5\textwidth]{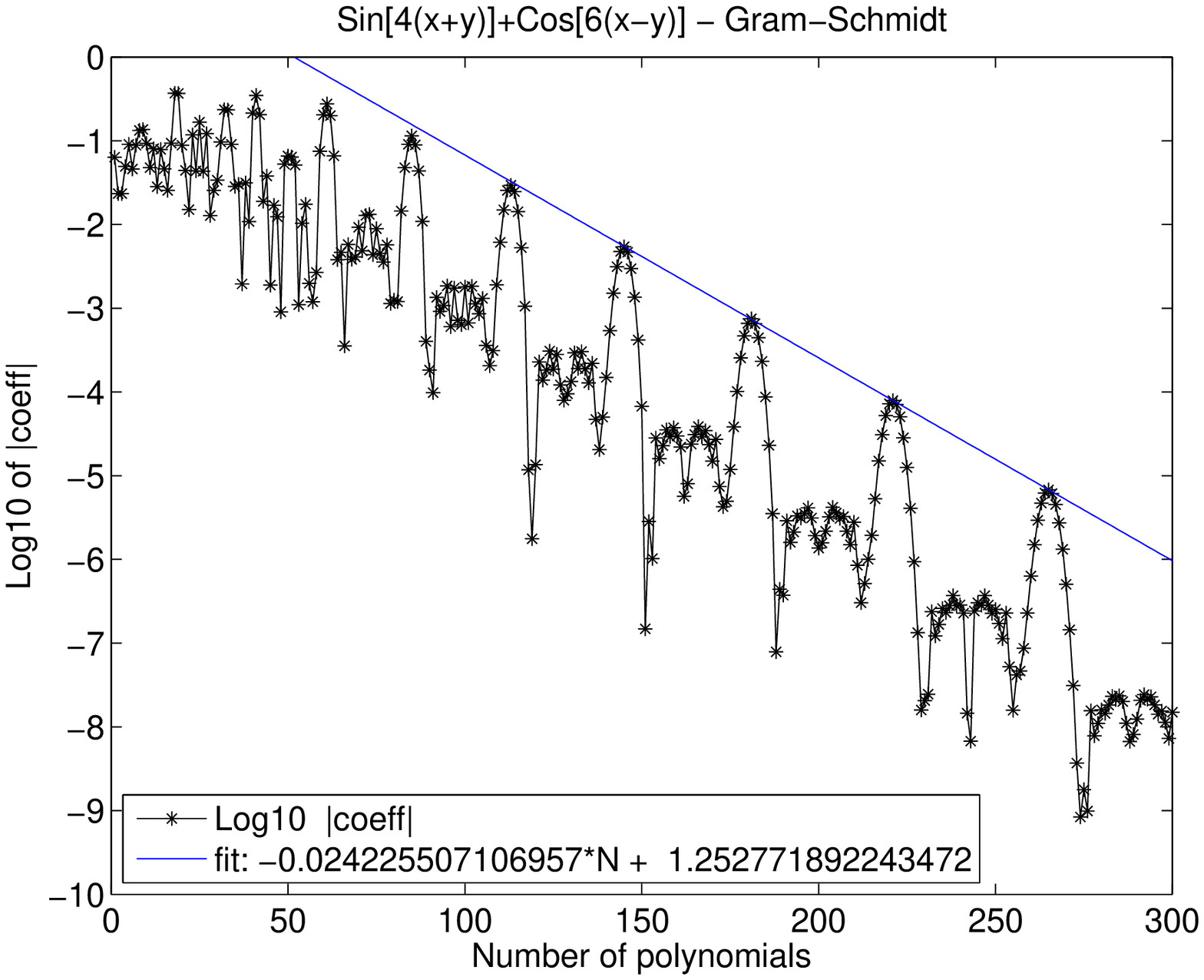}
\end{tabular}
\end{center}
\caption{Comparison between the coefficients of the expansions of $ f(x,y)=\sin\left(4(x+y)\right)+\cos(6(x-y))$ in a: Legendre polynomials and b: the orthogonal polynomials associated with$w_2(x,y)=({2}/{9}) \left(\chi_Q(x,y)+1\right)$.
} \label{fig:fig_num_exam_3.1}
\end{figure}

{\bf Example 3}: In this 2-D example we compare the rate of decay of the projection coefficients for the function $ f(x,y)=\sin\left(4(x+y)\right)+\cos(6(x-y))$.\newline \\
We compare the expansions of this function with respect to Legendre-type polynomials ($w_1(x,y)={1}/{4}$) and non-standard orthogonal polynomials ($w_2(x,y)=({2}/{9}) \left(\chi_Q(x,y)+1\right)$), where $\chi_Q(x,y)$ is the indicator function of the set.
\begin{equation}\label{Numerical_example_3_10} \nonumber
Q=\left \{  (x,y) \left | -\frac{1}{2}\leq x \leq\frac{1}{2}, \;-\frac{1}{2}\leq y \leq -x  \right . \right \} \;
\end{equation}
As before, he latter were obtained via Gram-Schmidt orthogonalization. The plot of the logarithm of the coefficients is presented in  Figure \ref{fig:fig_num_exam_3.1}. Note that although $w_2(x,y)$ is discontinuous, the decay rates of the coefficients are almost identical and that the peaks correspond to the same indexes, i.e. to polynomials which are generated by orthogonalization of the same monomials. The linear fits, generated by the last four peaks (polynomials 145, 181, 221 and 265) for $w_1$ and $w_2$ are $-0.024211589198446 N +  1.255763616975428$ and $-0.024225507106957 N +  1.252771892243472$, respectively.   

\begin{figure}[h]
\begin{center}
\begin{tabular}{ccc}
%a:Legendre&b:Gram Schmidt\\
\includegraphics[width=0.5\textwidth]{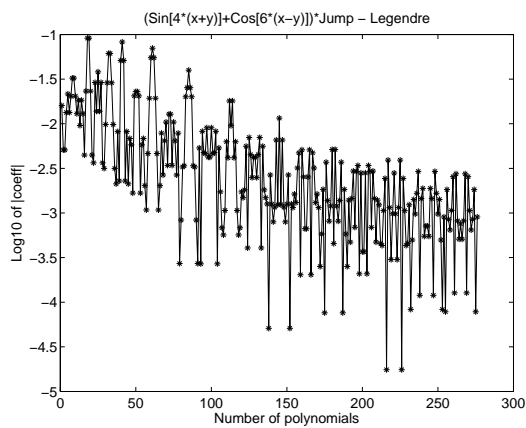}
\end{tabular}
\end{center}
\caption{Coefficients of the expansions of \newline $g(x,y)=[\sin(4(x+y))+\cos(6(x-y))]\cdot \frac{2}{9} \left(\chi_Q(x,y)+1\right)$ in Legendre type polynomials.
} \label{fig:fig_num_exam_3.20}
\end{figure}

We also calculate the coefficients of the expansion of 
\begin{equation}\label{Numerical_example_3_20} \nonumber
g(x,y)=[\sin(4(x+y))+\cos(6(x-y))]\cdot \frac{2}{9} \left(\chi_Q(x,y)+1\right)
\end{equation}
in 2-D Legendre-type polynomials ($w_1(x,y)=\frac{1}{4}$). Since $g$ is discontinuous the decay rate of the coefficients is much lower, see Figure \ref{fig:fig_num_exam_3.20}.
This example illustrates the potential use of nonstandard weight functions as a tool for obtaining fast coefficient decay rates for discontinuous or singular functions. In this case, instead of expanding, or integrating, $g$ directly one can use the discontinuous part, $\frac{2}{9} \left(\chi_Q(x,y)+1\right)$ as a weight function and obtain a much higher convergence rate, see Figure \ref{fig:fig_num_exam_3.1}b. This idea is not new. The use of expansions in Chebyshev, or Jacobi polynomials in order to get high convergence rates for integrating functions with singularities near the boundaries in the 1-D case was demonstrated e.g. in \cite{conte1980elementary}.  The lemma, however, provides the theoretical background for the  use of a larger class of weight functions in multidimensional settings. 

{\bf Example 4 - L-shaped domain}: We utilize the proposed approach for approximating the function
\begin{equation}
f(x,y)=\cos(x+y)\cdot(x^2+y^2)^\frac{1}{4}
\end{equation}

on the L-shaped domain $\Omega = \Omega_1 \cup \Omega_2$, where
\begin{eqnarray}
\Omega_1 &=& \left\{(x,y)|-1<x<1 , -1<y<0\right\} \\
\Omega_2 &=& \left\{(x,y)| 0<x<1 , -1<y<1\right\} \nonumber
\end{eqnarray}

With respect to the standard $L^2(\Omega)$ norm. Here, the approximated function can be written as \begin{equation}
f(x,y)=f_R(x,y)\cdot f_S(x,y)
\end{equation}
where $f_R(x,y)=\cos(x+y)$ is a regular function and $f_S(x,y)=(x^2+y^2)^\frac{1}{4}=w(x,y)$ can be regarded as a bounded and positive weight function on $\Omega$. Moreover, the partial derivatives of $f_S$ are unbounded near the origin. Thus, the origin is a point of weak singularity of $f$. Similar models frequently appear in heat transfer or elasticity problems, e.g \cite{szabo1996numerical}.

The approximation was done by decomposing $\Omega$ into a finite collection of equal squares $\left\{\Omega_j\right\}_{j=1}^{27}$ and constructing local approximations, $H_j:\Omega_j\rightarrow \mathbb{R}$, which define a global approximating function, $H$, by
\begin{equation}
H(x,y)=H_j(x,y) \, ; \, (x,y)\in \Omega_j
\end{equation}
We compare the following two approaches, for a given choice of the degree of the local approximation, $N_j$:
\begin{enumerate}
\item The 'classical' approach: the approximating functions, $H_j$, is constructed by using classical Legendre-type polynomials on $\Omega_j$ with respect to the inner product
\begin{equation}
\left<g,l\right>_{v_j}:= \iint_{\Omega_j} g(x,y)l(x,y)v_j(x,y)\,dx dy
\end{equation}
\begin{equation}
v_j(x,y):=\frac{1}{d_j} \, \,; \, \, d_j:=\mu(\Omega_j) 
\end{equation}
where $\mu$ denotes the standard Lebesgue measure on $\mathbb{R}^2$. The error on the square, $E_j:=\left|\left|f-H_j\right|\right|_{L^2(\Omega_j)}$, can be bounded by using convergence results for classical orthogonal polynomials. The total $L^2(\Omega)$ error is calculated by summation over the local errors.

\item The modified approach: the approximation on $\Omega_j$ is carried out with respect to the following inner product on $\Omega_j$ 
\begin{equation}
\left<g,l\right>_{w_j}:= \iint_{\Omega_j} g(x,y)l(x,y)w_j(x,y)\,dxdy
\end{equation} 
\begin{equation}
w_j(x,y):=\frac{w(x,y)}{d_j} \, \,; \, \, d_j:=\iint_{\Omega_j} w(x,y)\,dx = \iint_{\Omega_j} f_S(x,y) \,dxdy
\end{equation}
where $f_S$ plays the role of a weight function. Let $P_{N_j}f_R$ be the polynomial $L^2_{w_j}(\Omega_j)$-approximation of $f_R$ with respect to the polynomial basis that is generated by orthogonalization, $\left\{\Psi_k\right\}_{|k|\leq N_j}$. Thus, we have the following representation
\begin{equation}
P_{N_j}f_R=\sum_{|k|\leq N_j} \hat{f_R}_k \Psi_k \, \,; \, \, \hat{f_R}_k=\left<f_R,\Psi_k\right>_{w_j}.
\end{equation}
Hence, we define
\begin{equation}
H_j(x,y):=P_{N_j}f_R(x,y)\cdot f_S(x,y) \, ; \, (x,y)\in \Omega_j
\end{equation}
The local error on $\Omega_j$ is $E_j:=\left|\left|f-H_j\right|\right|_{L^2(\Omega_j)}$ and the total error is obtained by summation over the squares.
\end{enumerate}

\textbf{Remark:} Generating non-standard orthogonal polynomials for a non-classical weight function is a costly procedure. Thus, a more advanced approach is to use the approximations of the second type only for elements that are adjacent to the singularity.

We postpone the formal error analysis of the latter approach for future work which will also include a comparison to the full h-p approximation method. In this example, we used local approximations by cubic polynomials. 

\begin{figure}[h]
\begin{center}
\begin{tabular}{ccccclllll}
%a: Integration error&b: Collocation points\\
\includegraphics[width=1\textwidth]{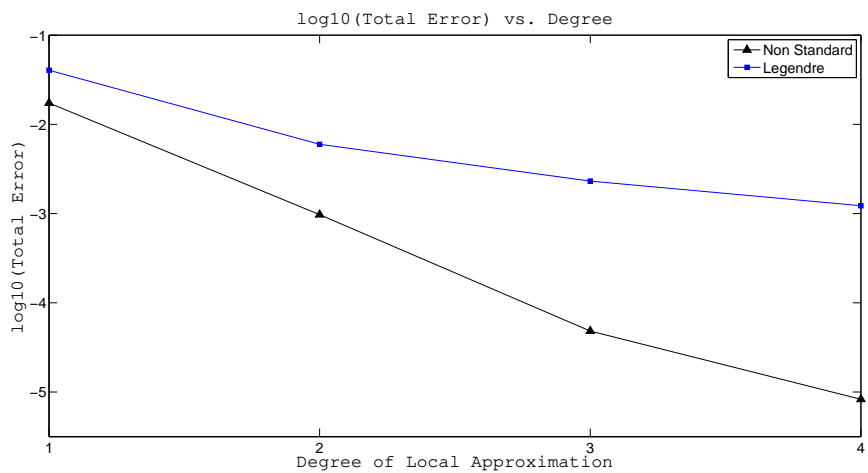}
\end{tabular}
\end{center}
\caption{L-Shaped - Example 3 - Total Error.
} \label{fig:LShex3E}
\end{figure}

\begin{figure}[h]
\begin{center}
\begin{tabular}{ccccclllll}
%a: Integration error&b: Collocation points\\
\includegraphics[width=1.1\textwidth]{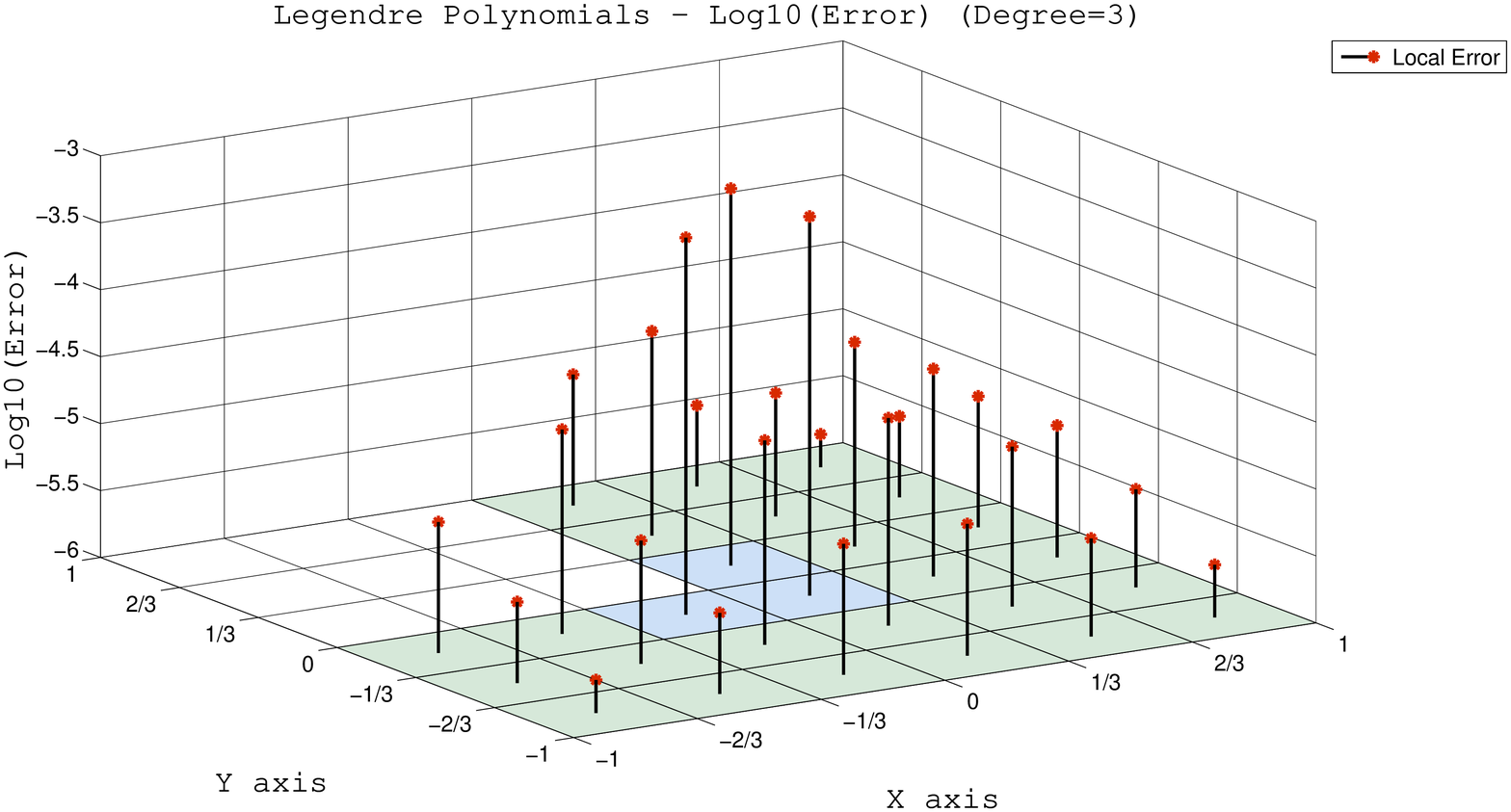}
\end{tabular}
\end{center}
\caption{L-Shaped - Example 3 - Legendre Polynomials , Local Error.
} \label{fig:LShex3B}
\end{figure}

\begin{figure}[h]
\begin{center}
\begin{tabular}{ccccclllll}
%a: Integration error&b: Collocation points\\
\includegraphics[width=1.1\textwidth]{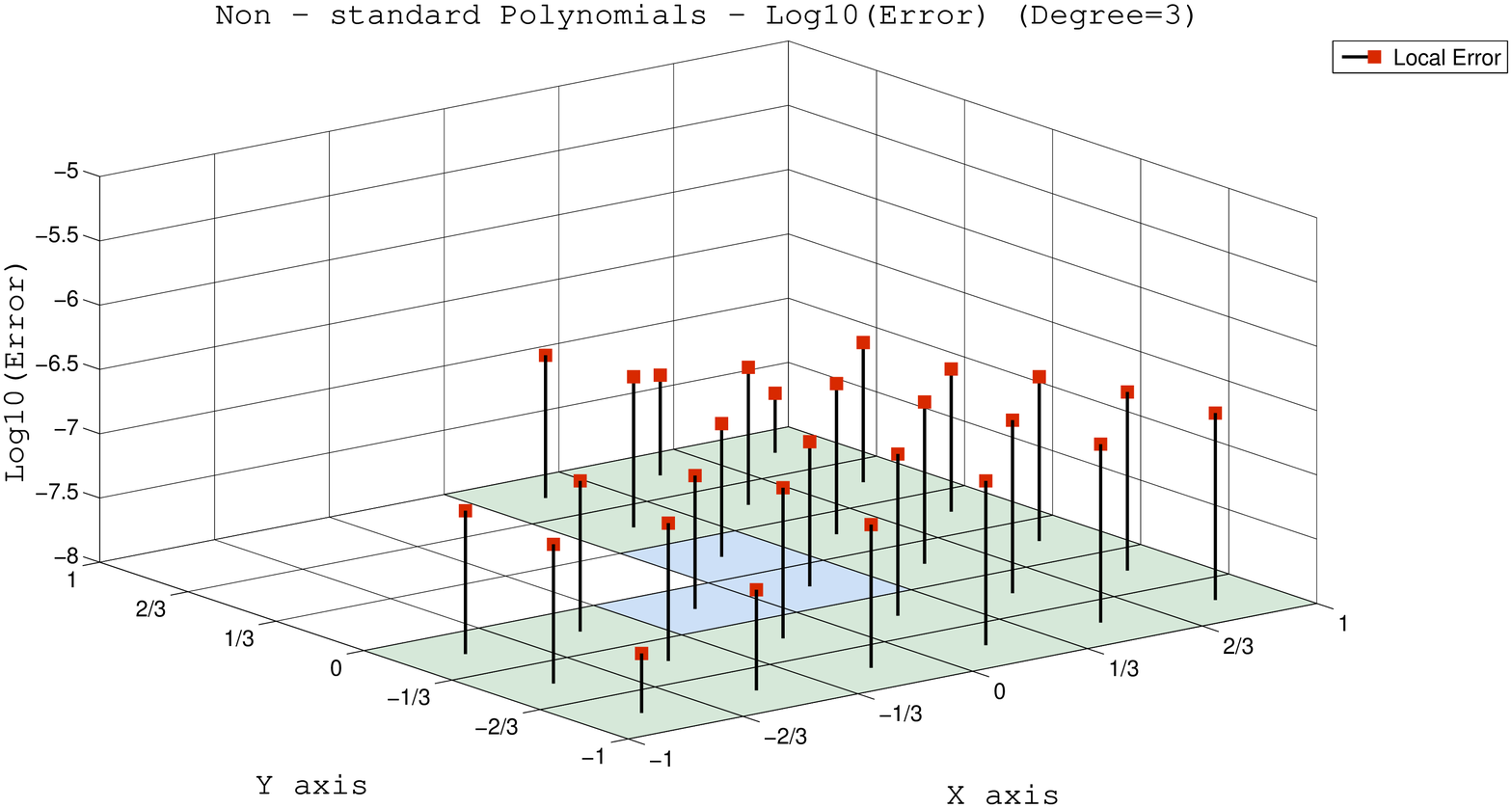}
\end{tabular}
\end{center}
\caption{L-Shaped - Example 3 - Non-standard Polynomials , Local Error.
} \label{fig:LShex3C}
\end{figure}

\begin{figure}[h]
\begin{center}
\begin{tabular}{ccccclllll}
%a: Integration error&b: Collocation points\\
\includegraphics[width=1.1\textwidth]{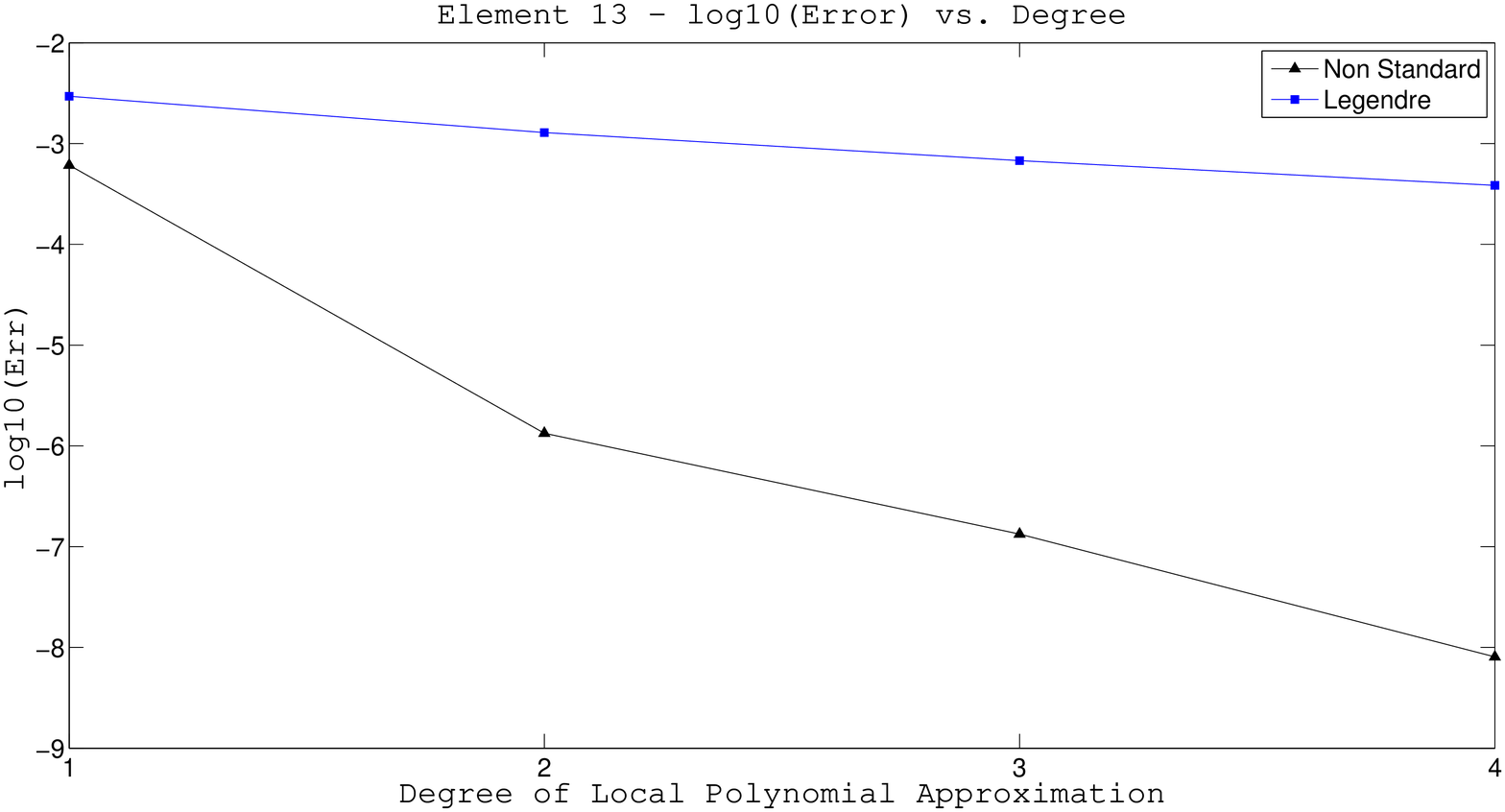}
\end{tabular}
\end{center}
\caption{L-Shaped - Example 3 - Element 13 , Local Error.
} \label{fig:LShex3FF}
\end{figure}

The domain decomposition and the numbering of the subdomains $\Omega_j$ are shown in figure (\ref{fig:LShex3A}). The errors on each element for local approximations of various degree are summarized in Appendix A. The total error over all elements is presented in figure (\ref{fig:LShex3E}).
\begin{figure}[h]
\begin{center}
\begin{tabular}{ccccclllll}
%a: Integration error&b: Collocation points\\
\includegraphics[width=1\textwidth]{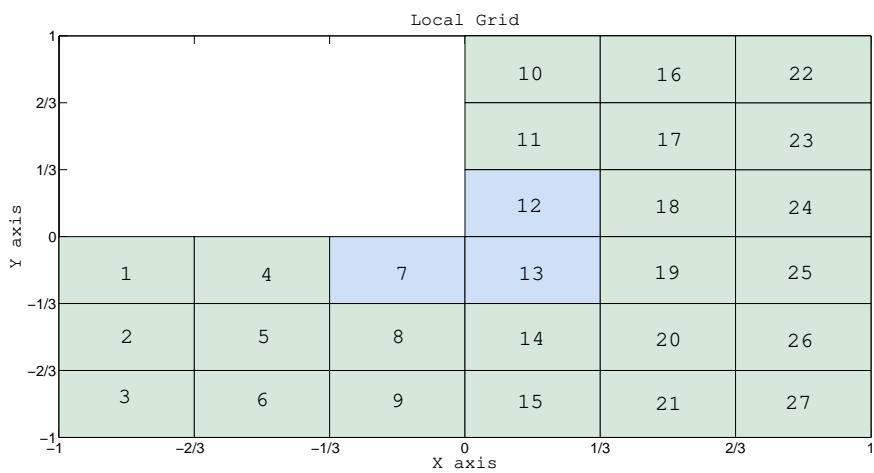}
\end{tabular}
\end{center}
\caption{L-Shaped - Example 3 - Domain Decomposition and the numbering of the subdomains $\Omega_j$.
} \label{fig:LShex3A}
\end{figure}
The local error for approximation on each element is shown in (\ref{fig:LShex3B}) and (\ref{fig:LShex3C}).
We also present the approximation error for element number 13, which is adjacent to the singularity, see figure (\ref{fig:LShex3FF}).
\bigskip \newline
These results show that the singularity of $f$ at the origin causes 
the approximation
errors in the Legendre-type polynomials expansion to be approximately 2.5 orders of magnitude larger than the errors in expansion using non-standard polynomials.  
%%%%%%%%%%%%%%%%%%%%%%%%%%%%%%%%%%%%%%%%%%%%%%%%%%%%%%%%%%%%%%%%%%%%%%%%%%%%%%
\FloatBarrier
\bigskip

%%%%%%%%%%%%%%%%%%%%%%%%%%%%%%%%%%%%%%%%%%

\section{Numerical Integration - A Heuristic Approach}\label{section:numerical integration}
Since it is rarely possible to symbolically integrate arbitrary functions, it is essential to provide an approach for constructing numerical integration formulas. Unfortunately, unlike the 1-D case, where Gauss quadrature can be used, there is no systematic way of generating such formulas in a multidimensional setting.
In this section, we demonstrate a possible approach to finding suitable numerical integration formulas for spectral integration.

The proposed approach is based on selection of collocation points that minimize the condition number of the induced cubature formula. This section presents the theoretical framework for this construction, followed by several numerical examples which demonstrate the validity of the described approach. The development of more efficient and robust methods of implementing the presented ideas is postponed for future work. 

Let $\Omega \in \RR^d$ and $w(x)$ be a weight function on $\Omega$, $f:\Omega \longrightarrow \RR$ and $f \in H_w^p(\Omega)$ for $p \ge 1$  where $ H_w^p(\Omega)$ is a weighted Sobolev space defined on $\Omega$.

\bigskip

We use the following notations:
\begin{itemize}
\item $\{\Psi_n(x)\}_{n \geq 0}$\, - \,Orthogonal polynomials on $\Omega$ with respect to the weighted inner product $\int_{\Omega} f(x)g(x)w(x) dx$. 
\item $P_N^{\Psi},Q_N^{\Psi}$\, - \,The projection operators that correspond to  $\{\Psi_n(x)\}_{n \geq 0}$.
 \item Let $ x_0,  \ldots, \,x_M$ be a set of points  in $\Omega$. We denote by $\vf$ the vector 
    \begin{equation} \label{eq_3.10}
    \vf\, = \,\begin{pmatrix} f(x_0)&\\ \vdots\\f(x_M) \end{pmatrix} \in \mathbb{R}^{M+1} 
    \end{equation}
\end{itemize}

{\bf Remark}: In this section, and the rest of the manuscript, $M+1$ is the number of polynomials in the constructed basis of the space of polynomials of degree less than or equal to $N$. In $d$ dimensions
    \begin{equation} \label{eq_3.12}
    (M+1)\, = \, {N+d \choose d }\, = \, O(N^d) \;.
    \end{equation}

In order to simplify the notations, we make a technical assumption that 
\begin{equation}\label{w_normalization}
 \int_{\Omega}  w(x) dx \, =\, 1  \;,
\end{equation}
otherwise, it is always possible to normalize $w$.

We assume that
\begin{equation}\label{Xiu_estimate} 
\left\|Q_N^{\Psi}f \right\|_{L_w^{2}}\, \leq \,const \cdot N^{-p} \cdot \left\|f\right\|_{H_w^{p}[\Omega]}
\end{equation}

In the case of non-standard orthogonal polynomials, this is guaranteed by lemma (\ref{main_theorem}). Let $ x_0,  \ldots, \,x_M \in\Omega$ be a set of collocation points. We would like to  find a numerical integration formula of the form:
\begin{equation}\label{eq_3.30}
I[g] := \int_{\Omega} g(x)w(x) d\Omega\, \approx \, \sum_{j=0}^{M} A_jg(x_j)\, = \, \left<\vA,\vg\right>  \;,
\end{equation}
where
\begin{equation} \label{eq_3.40}
\vA\, = \, \begin{pmatrix} A_0\\ \vdots\\ A_M \end{pmatrix} \;.
\end{equation}
We denote this numerical integration operator by $NI[g]$=$NI_{M+1}^{x_0\, ... \,x_M}[g]$.

If the collocation points $ x_0,  \ldots, \,x_M \in\Omega$ are known (for example, the Gauss or Gauss-Lobato collocation points) then the coefficients vector, $\vA$, could be found by demanding that $NI[g]$ is exact for $\{\Psi_0,...,\Psi_M\}$. This gives rise to the following system of linear equations
\begin{equation}\label{algebraic_precision}
    \begin{pmatrix} \Psi_0\left(x_0\right)& \Psi_0\left(x_1\right)& \cdots & \Psi_0\left(x_M\right) \\
    \Psi_1\left(x_0\right)& \Psi_1\left(x_1\right)& \cdots & \Psi_1\left(x_M\right)\\
    \vdots& \vdots& \cdots& \vdots\\
    \Psi_M\left(x_0\right)& \Psi_M\left(x_1\right)& \cdots & \Psi_M\left(x_M\right)
    \end{pmatrix}  \vA\, = \, \begin{pmatrix} \beta_0 \\
    0 \\
    \vdots \\
    0 \end{pmatrix}\,  \; ,
    %; \, \beta_0\in\mathbb{R}
\end{equation}
where
\begin{equation}\label{def_of_beta_0}
\beta_0 \, = \, \int_{\Omega} 1  \, w(x) dx  \, = \, 1 \,,
\end{equation}

which can be written as:
\begin{equation}\label{eq_3.60}
R \vA \, = \, \vbeta  , \;\;\;\text{or}  \;\;\;\vA=R^{-1}\vbeta
\end{equation}
and $R^{-1} = R^{-1}(x_0,\ldots,x_M)$.

The numerical integration formula is then given by
\begin{equation}\label{numerical_integration_formula} 
NI[g]\, = \, \left<\vA,\vg\right> \,=\, \left<R^{-1}\vbeta, \, \vg\right>
\end{equation}
Following \cite{quarteroni2010numerical} we define the absolute condition number of the numerical integration formula $NI$, $K_{ abs}$, by
\begin{eqnarray} \label{abs_condition_number_def}
K_{abs}(NI)
&:=& \,\sup_{\left\| \vu\right\|_2\,\leq \,1}\left|\left<R^{-1}(x_0,...,x_M)\vbeta,\vu\right>\right|  \nonumber \\
&\leq& \,\left\|R^{-1}(x_0,...,x_M)\vbeta\right\|_2 
\end{eqnarray}

\underline{\textbf{The proposed approach}}\bigskip \newline
For every $M\in \mathbb{N}$ find a set of collocation points in $\Omega$, denoted by $ x_0 ,\ldots ,x_M $, such that the bound $\left\|R^{-1}(x_0,...,x_M)\vbeta\right\|_2$ on $K_{abs}(NI)$ is as small as possible. We denote is bound by $\lambda(M)$, i.e
\begin{equation}\label{main_approach_numerical_integration} 
% \left\|NI_{M+1}^{x_0\, ... \,x_M}\right\|_2\, \leq \,
\left\|R^{-1}(x_0,...,x_M)\vbeta\right\|_2\, \leq \lambda(M)
\end{equation}
This is done by using a suitable minimization algorithm on $\left\|R^{-1}(x_0,...,x_M)\vbeta\right\|_2$ as a function of $M+1$ variables.\bigskip \newline
\textbf{Remarks:}
\begin{itemize}
\item We do not claim that an absolute minimum for $\left\|R^{-1}(x_0,...,x_M)\vbeta\right\|_2$ is required for the approach to work. Subsequent analysis will show that as long as $\lambda(M)$ does not grow too rapidly as $M$ tends to infinity, this approach will result in spectral convergence  of the resulting cubature formulas. Specifically, estimate (\ref{error_estimate_final_1}) shows that if $\lambda(M) = O\left (M^{\alpha} \right)$ where $\alpha$ is small with respect to the smoothness of the function, $p$, then the use of numerical integration, instead of analytic one will allow for convergence, although at at a slower rate.
\item It is well known that in the multidimensional case, system (\ref{algebraic_precision}) doesn't have to be invertible. However, for such collocation points, the condition number is infinite. Therefore, such points are avoided according to our approach (as are points for which the condition number is very large).
\item The collocation points can be found in a sequential manner, where in each step the previously found points are fixed and the minimization is carried out with respect to the new (additional) point.
\end{itemize}

\subsection{The construction method}
The proposed heuristic approach can be summarized as follows:

Given the Orthogonal Polynomials $\left\{\Psi_j\right\}_{j \ge 0}$:
    \begin{enumerate}
        \item Normalize $w$.
        \item Construct the matrix $R$ as a function of $x_0,...,x_M$ and find $x_0,...,x_M$ which minimize $\left\|R^{-1}(x_0,...,x_M)\vbeta\right\|_2$. This is done numerically using  an optimization algorithm.
        \item The obtained $x_0,...,x_M$ are the  integration points, $\vA=R^{-1}\vbeta$ are the integration weights and $NI[g]\, = \, \left<\vA,\vg\right> $, See \eqref{numerical_integration_formula}. 
    \end{enumerate}
In simple cases $\left\{\Psi_j\right\}_{j \ge 0}$ can be found analytically.  Otherwise, they can be computed using Gram Schmidt procedure or via a recursion formula that is based on Favard's theorem, using classical numerical integration schemes. Although this is a costly procedure it is done once as a preprocessing stage.

\subsection{Analysis}
The error of the cubature formula is defined by:
\begin{equation}\label{error_term}
     E\left[f\right] \, := \, I\left[f\right]-NI\left[f\right]\; ,\;\;\;\;\; \,f\in H_w^{p}[\Omega] \;.
\end{equation}
Since  $f=P_M^{\Psi}f+Q_M^{\Psi}f$ and  $NI\left[\cdot \right]$ is exact for the polynomials $\left\{\Psi_j\right\}_{j \ge 0}^M$
\begin{eqnarray}\label{error_term2}
\left|E\left[f\right]\right|\, &=& \,\left|I\left[P_M^{\Psi}f\right]+I\left[Q_M^{\Psi}f\right] \right .  \nonumber \\
&& \hspace{1cm}\left .-NI\left[P_M^{\Psi}f\right]-NI\left[Q_M^{\Psi}f\right]\right|\nonumber\\
&=&\left|I\left[Q_M^{\Psi}f\right]-NI\left[Q_M^{\Psi}f\right]\right|\nonumber\\
&\leq& \left|I\left[Q_M^{\Psi}f\right]\right|+\left| NI\left[Q_M^{\Psi}f\right]\right|
\end{eqnarray}

Using assumption \eqref{Xiu_estimate} we have that.
\begin{eqnarray}\label{symbolic_projection}
\left|I\left[Q_M^{\Psi}f\right]\right|&=&\left|\int_{\Omega} \left(Q_M^{\Psi}f\right)(x)w(x)dx\right|\nonumber\\
&=&\int_{\Omega}\left|\left(Q_M^{\Psi}f\right)(x)\right|\sqrt{w(x)}\sqrt{w(x)}dx \nonumber \\
&\leq& \sqrt{\int_{\Omega} \left|\left(Q_M^{\Psi}f\right)(x)\right|^2w(x)dx}\cdot \sqrt{\int_{\Omega}w(x)dx}\nonumber \\
&=&\left\|Q_M^{\Psi}f\right\|_{L_w^2} \, \leq \,const \cdot N^{-p} \cdot \left\|f\right\|_{H_w^{p}(\Omega)} 
\end{eqnarray}
Therefore $\left|I\left[Q_M^{\Psi}f\right]\right|$ tends to zero at a rate of $N^{-p}$ as $N$ tends to infinity. 

Evaluating the second term in \eqref{error_term2},  $\left| NI\left[Q_M^{\Psi}f\right]\right|$,  is much more difficult as it requires to estimate $\left| Q_M^{\Psi}f  \right|$ at the nodes  $ x_0,  \ldots, \,x_M$ rather than it's integral norm. It is well know that the sequence $\left \| \Psi_j \right \|_{\infty}$ can be unbounded as $j \longrightarrow \infty$. For example, the $L_2$ normalized Legendre polynomials satisfy   $\left \| P_j \right \|_{\infty} = O(\sqrt{j})$. For general orthogonal polynomials there are no known evaluations of  $\left \| \Psi_j \right \|_{\infty}$. Another approach could be the use of inequalities which bound the maximum norm of polynomials, or some norm of their derivatives, by the norm of these polynomials. Well-known inequalities are the  Bernstein and Nikolskii inequalities, see e.g. \cite{devore1993constructive}. However, they are not directly applicable in this setting due to the fact that there are no generalizations of these inequalities to several dimensions and to nonstandard weight functions. Therefore, we use the following approach:
\begin{itemize}
\item We first show that the Lebesgue measure (denoted by $\mu$) of the set in which  $\left| \left(Q_M^{\Psi}f\right)(x)\right| $ is larger than $M^{-\theta}$, for an auxiliary parameter $\theta>0$,  tends to zero as $M\longrightarrow \infty$.
\item We then show that when the nodes $ x_0,  \ldots, \,x_M$ are taken in the complement of this set the integration error is small. In practice, once $x_0,...,x_M$ are chosen, it can be verified that for a given $\theta$
\begin{equation}\label{error_estimate_final_verification}
\left|\left(Q_M^{\Psi}f\right)(x_j)\right| \, = \, \left|f(x_j) - \left(P_M^{\Psi}f\right)(x_j)\right|  \le M^{-\theta}\; , \;\; \forall j=0,...,M
\end{equation}
by substitution.
\item In the second part of the following analysis, there is an implicit assumption that $\mu(\Omega)<\infty$. In practice, this assumption poses no loss of generality. Indeed, since $f \in L_{2,w}(\Omega)$ we know (see \cite{stein2016harmonic}) that for a ball $B=B(0,R)$ with a sufficiently large radius 
\begin{equation}
\int_{B^c} |f(x)|^2w(x) dx \leq \epsilon^2
\end{equation}
where $\epsilon$ is the machine precision. Thus, it is sufficient to approximate the integral of the function $f(x)$ on the ball $B$ without incurring any loss of accuracy. This follows from 
\begin{eqnarray}\label{remarkdomain}
\left|\int_{B^c} f(x)w(x)dx\right|
&=&\int_{B^c}\left|f(x)\right|\sqrt{w(x)}\sqrt{w(x)}dx \nonumber \\
&\leq& \sqrt{\int_{B^c} \left|f(x)\right|^2w(x)dx}\cdot \sqrt{\int_{B^c}w(x)dx}\leq \epsilon
\end{eqnarray}
\end{itemize}
As we've mentioned, this is a theoretical analysis of the approach and, thus, the auxiliary parameter $\theta$ does not appear in any of the calculations carried out in practice. Moreover, we assume that f is sufficiently smooth so that such a parameter $\theta$ can be obtained.\bigskip \newline 
Note that
\begin{eqnarray}\label{numerical_integration_of_projection}
\left|NI\left[Q_M^{\Psi}f\right]\right| &\leq& \left\|NI\right\|_2 \,  \left\|  \vQ_{\vM}^{\vPsi}\vf\right\|_2 \nonumber\\
&\leq& \lambda(M) \cdot \sqrt{M} \cdot \max_{\,0\leq j \leq M}\left|\left(Q_M^{\Psi}f\right)(x_j)\right| 
\end{eqnarray}

We now divide the analysis into two cases:

{\bf Case 1}: $w(x) \ge \delta>0$

Let $\theta >0$.  By using Markov's inequality we obtain
\begin{eqnarray}\label{tail_local_control}
\mu\left(\left\{x \in \Omega\, | \,\left|\left(Q_M^{\Psi}f\right)(x)\right| > M^{-\theta}\right\}\right)
&\leq&  M^{2\theta} \int_{\Omega} \left|\left(Q_M^{\Psi}f\right)(x)\right|^2dx \nonumber \\
&=&\frac{M^{2\theta}}{\delta}\int_{\Omega}\left|\left(Q_M^{\Psi}f\right)(x)\right|^2\delta dx \nonumber \\
&\leq&\frac{M^{2\theta}}{\delta}\left\|Q_M^{\Psi}f\right\|_{L_w^2}^2\nonumber \\
 & \leq & \frac{M^{2\theta}}{\delta}  \,const \cdot N^{-2p} \cdot \left\|f\right\|^2_{H_w^{p}(\Omega)} \nonumber\\
 & \approx & \frac{N^{2(\theta d-p)}}{\delta}  \,const  \cdot \left\|f\right\|^2_{H_w^{p}(\Omega)} 
\end{eqnarray}
Therefore if $\theta d-p<0$

\begin{equation}\label{tail_local_control_2}
\mu\left(\left\{x \in \Omega\, | \,\left|\left(Q_M^{\Psi}f\right)(x)\right| > M^{-\theta}\right\}\right)
\overset{}{\underset{N \rightarrow \infty}\longrightarrow} 0  
\end{equation}
Thus, by selecting the collocation points in the complement of this set (which has measure of at least $1-\left (( {N^{2(\theta d-p)}}/{\delta}) \cdot const  \cdot \left\|f\right\|^2_{H_w^{p}(\Omega)} / \mu\left( \Omega \right) \right)$)
\begin{equation}\label{numerical_integration_of_projection}
\left|NI\left[Q_M^{\Psi}f\right]\right| 
\leq \lambda \cdot M^{1/2-\theta} \, \approx \, \lambda \cdot N^{(1/2-\theta)d} \;,
\end{equation}
Finally, we obtain that
\begin{equation}\label{error_estimate_final_1}
\left|E\left[f\right]\right|\, \leq \,const \cdot N^{-p} \cdot \left\|f\right\|_{H_w^{p}[\Omega]} + \lambda \cdot N^{(1/2-\theta)d} \rightarrow 0
\end{equation}
As $N$ tends to infinity. This ensures that the numerical integration formula approximates the exact integration for sufficiently large values of $N$.

{\bf Case 2}:  If there is no positive constant $\delta$ that bounds $w(x)$ from below we use a different approach. For every $n\in \mathbb{N}$ define the following set
\begin{equation}
S_n=\left\{x\in \Omega\, | \,w(x)<\frac{1}{n}\right\} \;.
\end{equation}
The sequence $\left\{S_n\right\}_{n\geq1}$ is downward monotone with respect to inclusion. Since $\mu(S_1) \leq \mu(\Omega)< \infty$, by applying the downward monotone convergence theorem, we obtain that
\begin{equation}\label{downward_sequence}
0=\mu\left(\left\{x\in \Omega\, | \,w(x)=0\right\}\right)= \lim_{n\to\infty}\mu\left(S_n\right) \;.
\end{equation}
Where we use the fact that $w(x)$ is a weight function that does not vanish on positive measure sets. Therefore, there exists a number $N_0 \in \mathbb{N}$ such that
\begin{equation}
\mu\left(S_{N_0}\right) < \eta
\end{equation}
For $\eta >0$ small enough. Therefore,
\begin{equation}\label{pointwise_control_case_2_part_1}
\mu\left(\left\{x\in \Omega\, | \,\left|\left(Q_M^{\Psi}f\right)(x)\right|>M^{-\theta},w(x)<\frac{1}{N_0}\right\}\right) \leq \eta 
\end{equation}
Moreover,
\begin{eqnarray}\label{pointwise_control_case_2_part_2}
\left\{x\in  \Omega\, | \,\left|\left(Q_M^{\Psi}f\right)(x)\right|>M^{-\theta},w(x) \geq \frac{1}{N_0}\right\}  \hspace{-8 cm}  \nonumber\\
&\subseteq&\left\{x\in  \Omega\, | \,\left|\left(Q_M^{\Psi}f\right)(x)\right|\cdot\sqrt{w(x)}>M^{-\theta}\frac{1}{\sqrt{N_0}},w(x) \geq \frac{1}{N_0}\right\}  \\
&\subseteq&\left\{x\in  \Omega\, | \,\left|\left(Q_M^{\Psi}f\right)(x)\right|\cdot\sqrt{w(x)}>M^{-\theta}\frac{1}{\sqrt{N_0}}\right\} \nonumber\\
&\subseteq&\left\{x\in  \Omega\, | \,\left|\left(Q_M^{\Psi}f\right)(x)\right|^2 \cdot w(x)>M^{-2\theta}\frac{1}{N_0}\right\}\nonumber
\end{eqnarray}
Thus, by using Markov's inequality and (\ref{Xiu_estimate}) we get the following estimate
\begin{eqnarray}\label{pointwise_control_case_2_part_3}
\mu\left(\left\{x\in  \Omega\, | \,\left|\left(Q_M^{\Psi}f\right)(x)\right| > M^{-\theta},w(x)\geq\frac{1}{N_0}\right\}\right)
 \hspace{-4 cm}  \nonumber \\
&\leq& N_0 \cdot M^{2\theta} \cdot \left\|Q_M^{\Psi}f\right\|_{L_w^2}^2  \\
&\leq & const \cdot N_0 \left\|f\right\|_{H_w^{p}(\Omega)}^2 M^{2 \theta } N^{-2p}\nonumber\\
&\approx & const \cdot N_0 \left\| f\right\|_{H_w^{p}(\Omega)}^2 N^{2\left(\theta d - p \right)} \nonumber
\end{eqnarray}
Finally, by using (\ref{pointwise_control_case_2_part_1}) and (\ref{pointwise_control_case_2_part_3}) we have
\begin{eqnarray}\label{pointwise_control_case_2_part_4}
\mu\left(\left\{x \in  \Omega\, | \,\left|\left(Q_M^{\Psi}f\right)(x)\right| > M^{-\theta}\right\}\right)  \hspace{-3cm}\nonumber \\
 &\leq & const \cdot N_o \left\|f\right\|_{H_w^{p}\Omega}^2 M^{2\left(\theta - p \right)} + \eta \nonumber \\
 &\approx& const \cdot N_o \left\|f\right\|_{H_w^{p}\Omega}^2 N^{2\left(\theta d - p \right)} + \eta 
\end{eqnarray}
Which implies that, with high confidence,  the points $ \{x_0,\ldots , x_M \}$ that minimize the condition number, will also   belong to the set 
$$\left\{x \in  \Omega\, | \,\left|\left(Q_M^{\Psi}f\right)(x)\right| < M^{-\theta}\right\} \;.$$
 As was noted above, this can be verified by substitution, see \eqref{error_estimate_final_verification}. 
The rest of the analysis and the conclusions are identical to the previous case.

\subsection{Numerical examples}
In this section, we present three examples of numerical integration using the algorithm described above. These examples demonstrate the validity of the described method, under the mentioned assumptions. For simplicity, the optimization was carried out using sequential choice of points via the downhill simplex method, see e.g.    \cite{nelder1965simplex}.

\begin{figure}[h]
\begin{center}
\begin{tabular}{lllll}
a: Integration error&b: Collocation points\\
\includegraphics[width=0.5\textwidth]{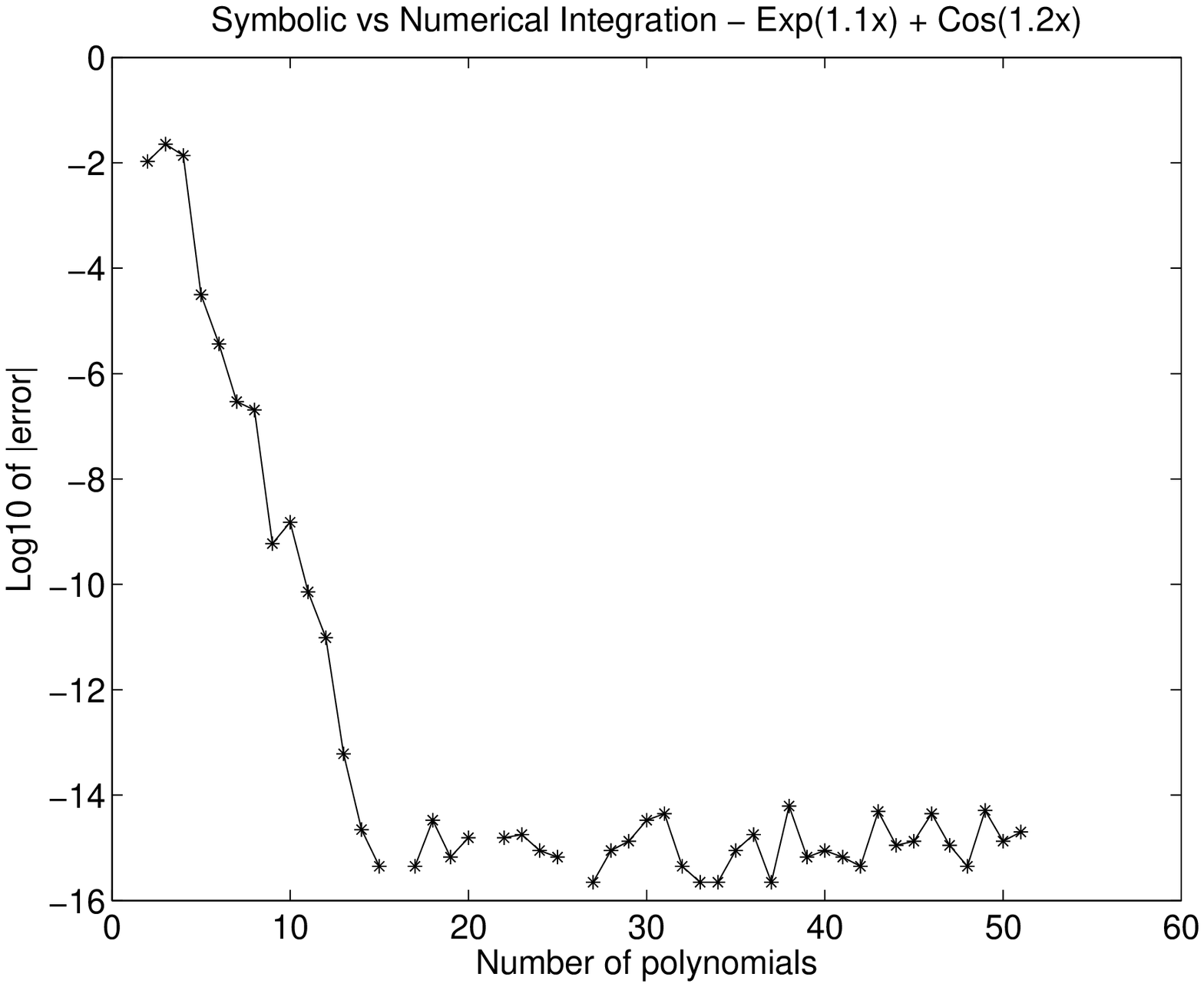}&
\includegraphics[width=0.55\textwidth]{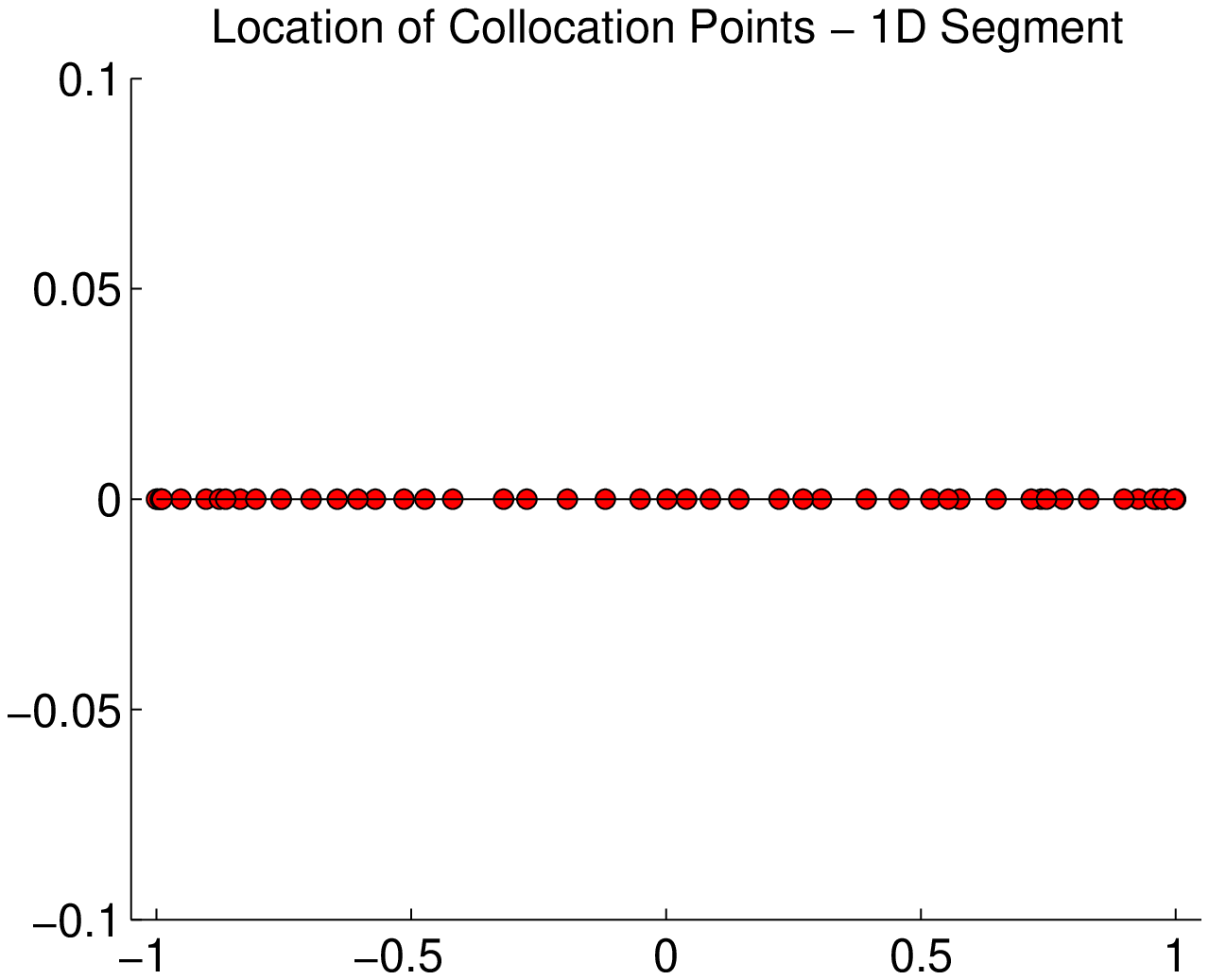}
\end{tabular}
\end{center}
\caption{1-D example. a: The integration error vs. the number of polynomials. b: the location of the collocation points.
} \label{fig:fig_num_exam_Integration_1.1}
\end{figure}

{\bf Example 1}: 
In this example we tested our approach to finding numerical integration formulas for the weight function  $w(x)=\frac{1}{3} \left(\chi_Q(x)+1\right)$, where $\chi_Q(x)$ is the indicator function of the set
 $Q=\{-\frac{1}{2}\leq x\leq\frac{1}{2} \}$
and the integration domain is $[-1,1]$. The integration was  carried out using the process presented above  for the function 
\begin{equation}
f(x)=e^{1.1x}+\cos(1.2x) \;.
\end{equation}
The $\log_{10}$ of the error vs the number of  collocation points is presented in Figure \ref{fig:fig_num_exam_Integration_1.1}a. The spectral convergence can be clearly seen from this plot. The locations of the collocation points are shown in Figure \ref{fig:fig_num_exam_Integration_1.1}b.  

\begin{figure}[h]
\begin{center}
\begin{tabular}{lllll}
a: Integration error&b: Collocation points\\
\includegraphics[width=0.5\textwidth]{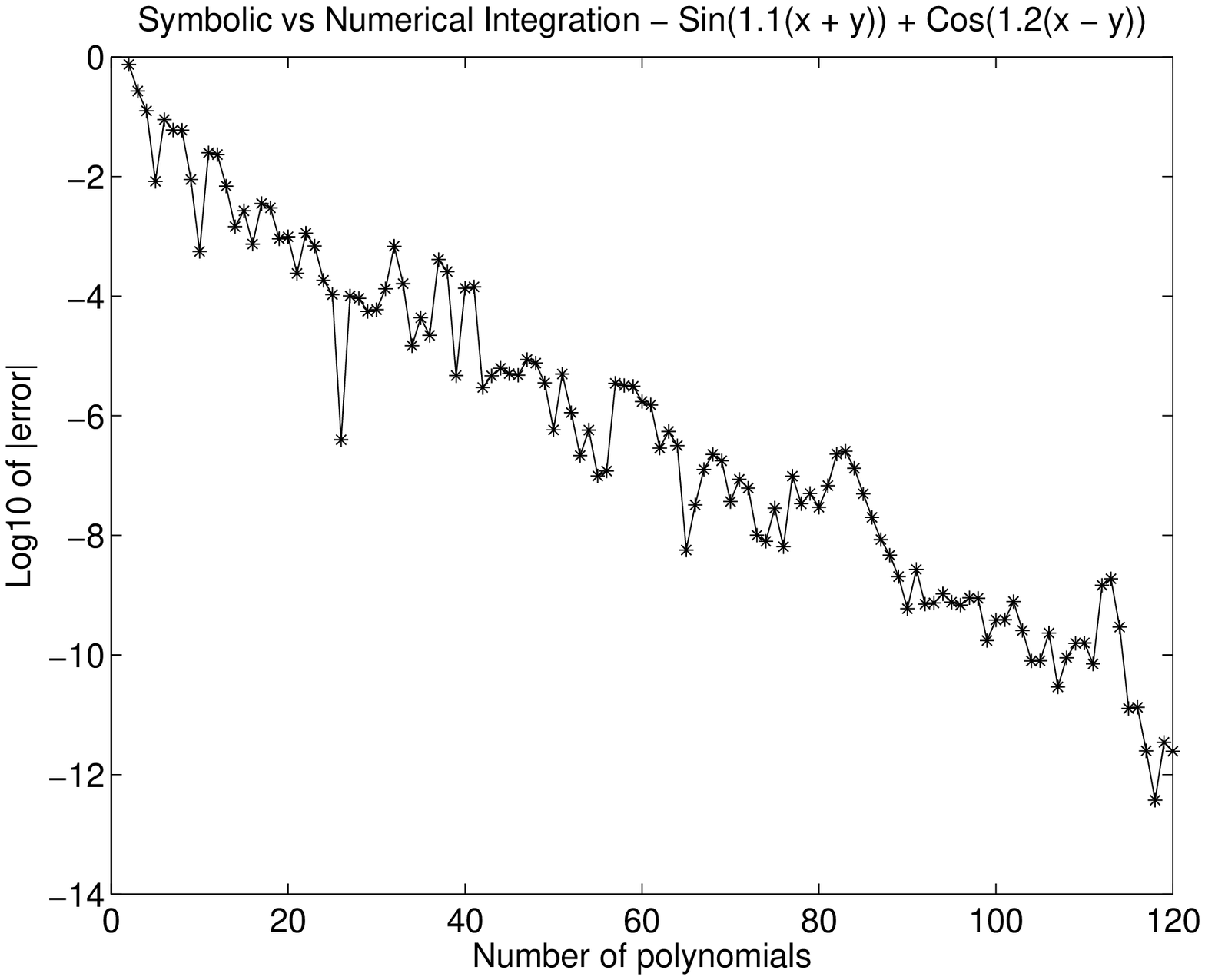}&
\includegraphics[width=0.55\textwidth]{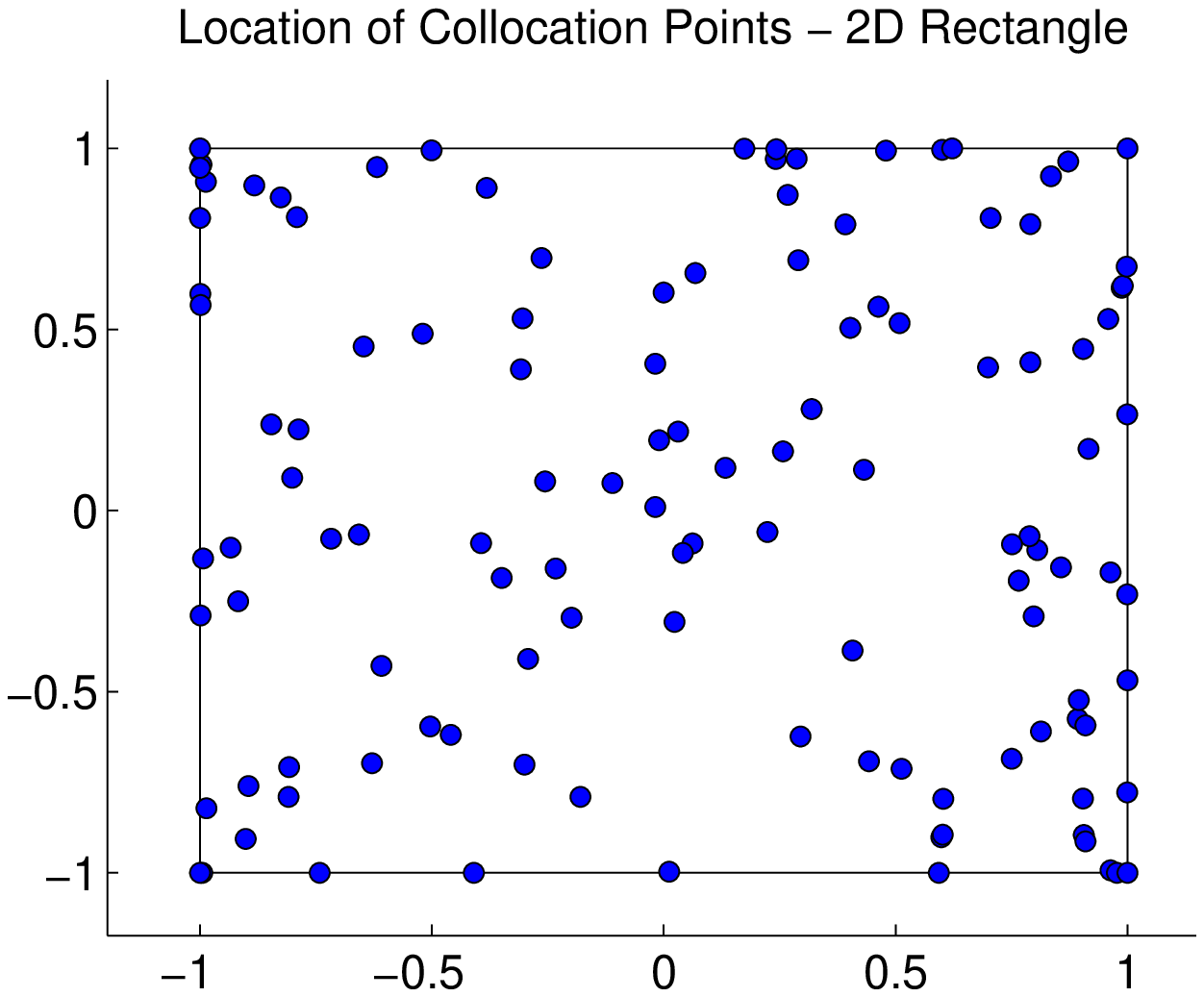}
\end{tabular}
\end{center}
\caption{2-D example. a: The integration error vs. the number of polynomials. b: the location of the collocation points.
} \label{fig:fig_num_exam_Integration_2.1}
\end{figure}

{\bf Example 2}: In this example  we integrate the function 
\begin{equation}
f(x,y)=\sin(1.1(x+y))+\cos(1.2(x-y))\;.
\end{equation} 
in the domain $(x,y)\in [-1,1]^2$ with the 
 weight function  $w(x,y)=\frac{2}{9} \left(\chi_Q(x,y)+1\right)$,  where $\chi_Q(x,y)$ is the indicator function of the set 
\begin{equation}
Q=\left\{-\frac{1}{2}\leq x\leq\frac{1}{2},-\frac{1}{2}\leq y \leq -x\right\} \;.
\end{equation}
The $\log_{10}$ of the error vs the number of  collocation points is presented in Figure \ref{fig:fig_num_exam_Integration_2.1}a. The spectral convergence can be clearly seen from this plot. The locations of the collocation points are shown in Figure  \ref{fig:fig_num_exam_Integration_2.1}b.

{\bf Example 3}: In this example  we integrate the function 
\begin{equation}
f(x,y)=\sin(1.1(x+y))+\cos(1.2(x-y))\;,
\end{equation} 
with the  weight function  $w(x,y)=2$ on the triangular domain 
\begin{equation}
\Omega=\left\{0 \leq x \leq 1,0\leq y \leq x\right\} \;.
\end{equation}

\begin{figure}[h]
\begin{center}
\begin{tabular}{lllll}
a: Integration error&b: Collocation points\\
\includegraphics[width=0.5\textwidth]{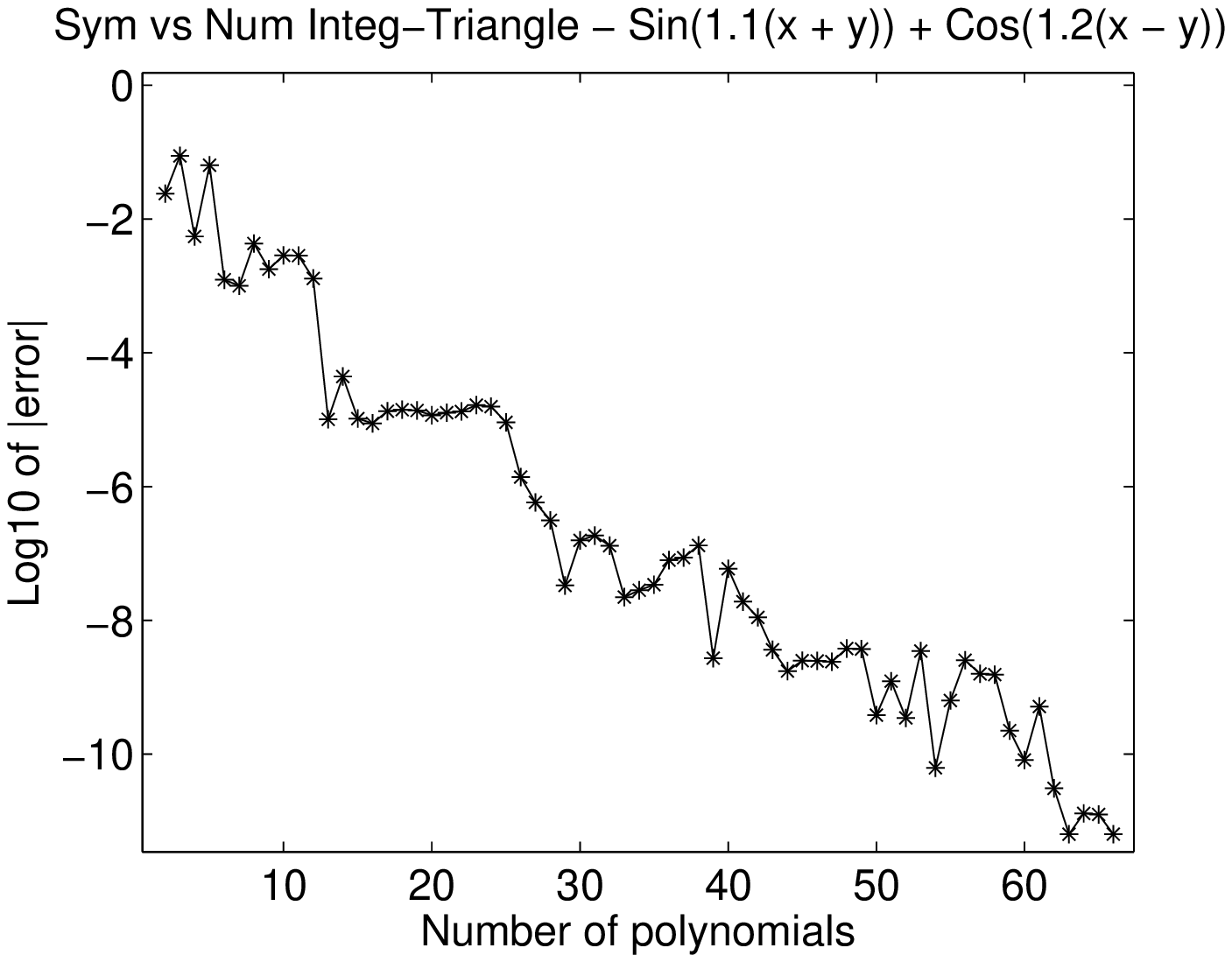}&
\includegraphics[width=0.55\textwidth]{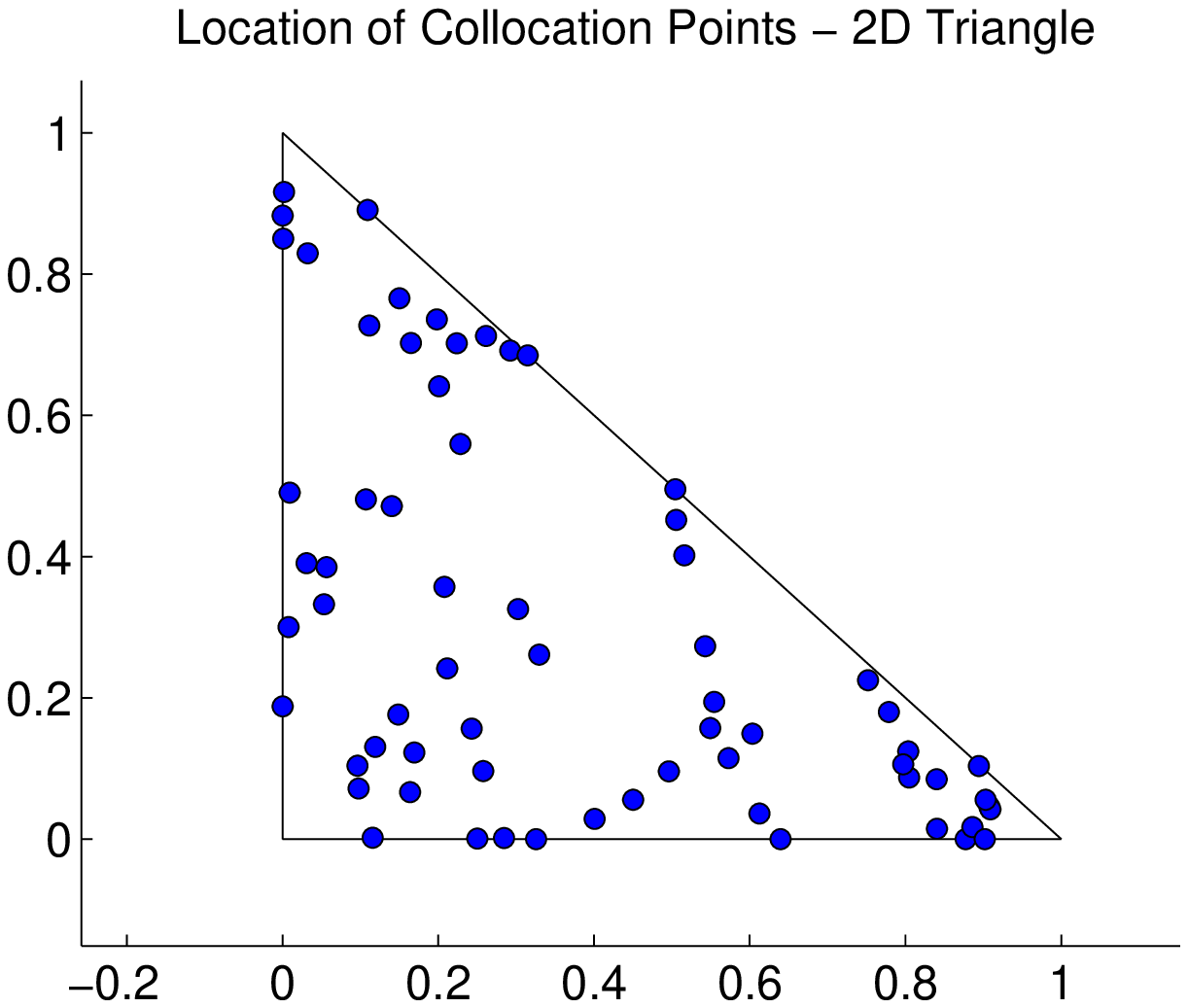}
\end{tabular}
\end{center}
\caption{2-D example. a: The integration error vs. the number of polynomials. b: the location of the collocation points.
} \label{fig:fig_num_exam_Integration_3.1}
\end{figure}

The $\log_{10}$ of the error vs the number of  collocation points is presented in Figure \ref{fig:fig_num_exam_Integration_3.1}a. As in the previous examples, the spectral convergence can be clearly seen from this plot. The locations of the collocation points are shown in Figure  \ref{fig:fig_num_exam_Integration_3.1}b.

These examples demonstrate the applicability of this integration method. It should be noted that though the locations of the collocation points look random, they are not. Choosing random points leads to a very large condition number and poor convergence rate.

\section{GPC Example}\label{Example_of_implementation_to_GPC}

In this section, we apply the theory presented above to solving a partial differential equation with stochastic parameters. We consider the following Cauchy problem, which is similar to Example 4.5 in \cite{xiu2010numerical}.
\begin{equation}\label{GPC_example_10}
u_t  =  a(z_1, z_2)^2u_{xx} \ \ \ t>0, x\in \mathbb{R} \ \ \  ;\ \ \  u(x,0)  =  e^{-x} \ \ \ x\in \mathbb{R} \;,
\end{equation}
where $(z_1, z_2)$ is a random vector with values in $[-1,1]^2$ and a probability density function
\begin{equation}\label{GPC_example_20}
w(z_1, z_2)=\frac{2}{9}(\chi_Q(z_1, z_2)+1) \ \ \ ; \ \ \  Q=\left\{-\frac{1}{2}\leq z_1\leq \frac{1}{2},-\frac{1}{2}\leq z_2 \leq -z_1\right\} 
\end{equation}
and $a(z_1, z_2)=\sqrt{\ln(\cos(z_1- z_2)+\sin(1.1(z_1,+z_2))+4)}$. The solution to this Cauchy problem is a function $u(x,t,z_1, z_2)$ which depends on the time and spatial variable, as well as on the random vector. Our goal is to estimate $u(0,1,z_1, z_2)$ by using stochastic collocation based on GPC. The solution to the initial value problem is
\begin{equation}\label{GPC_example_30}
u(0,1,z_1, z_2)= \cos(z_1-z_2)+\sin(1.1(z_1+ z_2))+4 \;.
\end{equation}

Denote by $\left\{\Phi_k(z_1, z_2)\right\}_{|k|\geq 0}$ the set of orthogonal polynomials on $[-1,1]^2$ with respect to $w(z_1, z_2)$, which were obtained via Gram-Schmidt algorithm. Thus, $u(0,1,z_1, z_2)$ can be expanded into a series given by
\begin{equation}\label{GPC_example_40}
u(0,1,z_1, z_2) = \sum_{|k|=0}^{\infty} \hat{u}_k(0,1)\Phi_k(z_1, z_2) 
\end{equation}
where
\begin{eqnarray}\label{GPC_example_50}
 \hat{u}_k(0,1) &=&\left(\Phi_k(\cdot,\cdot),\, u(0,1,\cdot,\cdot)\right)_{w}   \nonumber \\
 &=& \int_{-1}^{1}\int_{-1}^{1} \Phi_k(z_1, z_2)\, u(0,1,z_1, z_2) w(z_1, z_2) \, dz_1\, dz_2
\end{eqnarray}

For $N \in \mathbb{N}$ we define 
\begin{equation}\label{GPC_example_60}
\left(P_N \, u\right)(0,1,z_1, z_2)=\sum_{|k|=0}^{N} \hat{u}_k(0,1)\Phi_k(z_1, z_2)
\end{equation}
to be the projection of $u(0,1,z_1, z_2)$ on the set of orthonormal polynomials of degree less or equal to $N$. Instead of using the symbolic projections $\left(P_Nu\right)(0,1,z_1, z_2)$ we approximate $u(0,1,z_1, z_2)$ by numerical projections
\begin{equation}\label{GPC_example_70}
\left(\tilde{P}_Nu\right)(0,1,z_1, z_2)=\sum_{|k|=0}^{N} \tilde{u}_k(0,1)\Phi_k(z_1, z_2) \;,
\end{equation}
where the coefficients $\tilde{u}_k(0,1)$ are obtained by using a numerical integration formula based on $2N+2 \choose 2$   collocation points (the number of points is equal to the number of orthogonal polynomials of degree less or equal to $2N$). By the Pythagoras theorem, we have
\begin{eqnarray}\label{GPC_example_70}
 \left    \lVert u(0,1,\cdot,\cdot)-\left(\tilde{P}_Nu\right)(0,1,\cdot,\cdot) \right \rVert^2 &=& \left \lVert \left(\tilde{P}_Nu\right)(0,1,\cdot,\cdot)-\left(P_Nu\right)(0,1,\cdot,\cdot) \right \rVert^2 \hspace{-3 cm} \nonumber \\
    &+& \left \lVert u(0,1,\cdot,\cdot)-\left(P_Nu\right)(0,1,\cdot,\cdot) \right \rVert^2 \;
\end{eqnarray}
The decay rate of the second term can be estimated with the main lemma. Thus, it is sufficient to estimate the behavior of the first term for different values of $N \in \mathbb{N}$ . By the Parseval's equality, we have that the second term is equal to
\begin{equation}\label{GPC_example_80}
H(N):=\sum_{|k|=0}^{N} \left| \tilde{u}_k(0,1)- \hat{u}_k(0,1)\right|^2
\end{equation}

\begin{figure}[h]
\begin{center}
\begin{tabular}{ccccclllll}
%a: Integration error&b: Collocation points\\
\includegraphics[width=0.7\textwidth]{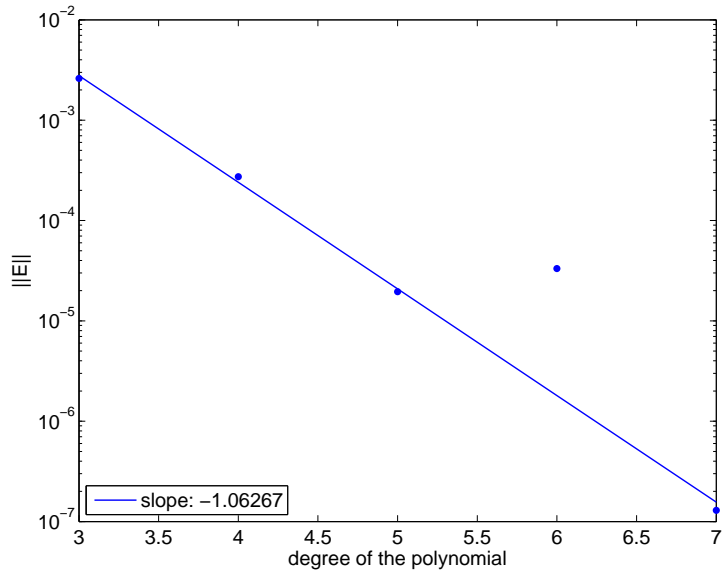}
\end{tabular}
\end{center}
\caption{plot of $\log_{10} H(N)$  vs $N$.
} \label{fig:GPC_convegrence_1}
\end{figure}

The  plot of $\log_{10} H(N)$  vs $N$ is presented in Figure  \ref{fig:GPC_convegrence_1}. As can be seen, $H(N)$ decays exponentially with $N$. The outlier in $N=6$ is probably due to the
unsophisticated optimization algorithm used for finding the collocation points. See further discussion in the conclusion section.

%\begin{center}
%\begin{tabular}{||c | c||}
%\hline
%N & $log_{10} H(N)$ \\
%\hline\hline
%2 & -1.934671 \\
%\hline
%3 & -2.583989 \\
%\hline
%4 & -3.563324 \\
%\hline
%5 & -4.709338 \\
%\hline
%6 & -4.479247 \\
%\hline
%7 & -6.887757 \\
%\hline
%\end{tabular}
%\end{center}

Thus, it can be seen that the error which arises from using the numerical integration formula instead of analytic one also decays. The analytic decay rate is faster than the symbolic approximation error $\lVert u(0,1,\cdot,\cdot)-\left(P_Nu\right)(0,1,\cdot,\cdot) \rVert^2$, as can be seen from the graph.

% \section{Conclusuins} \label{section:conclusuins}

%%%%%%%%%%%%%%%%%%%%%%%%%%%%%%%%%%%%%%%%%%%%%%%%%%%%%%%%%%%%%%%%%%%%%%%%%%

\section{Summary}\label{section:conclusions}
In this work, we showed that the convergence rate of expansions in orthogonal polynomials with  respect to a 'general' weight  function $w_2$ can be bounded by the convergence rate of  expansions in polynomials  associated with  respect to  weight  function $w_1$, as long as  $w_2$ is dominated by $w_1$. This result allows the evaluation of  the convergence rate  even in the case where the orthogonal polynomials are not eigenfunctions of some singular Sturm-Liouville
problem. This result gives the theoretical background for the use of these expansions. We also demonstrated that the use of nonstandard weight functions enables to derive fast convergent expansions of non-smooth functions.

We also presented an approach to derive a multi-dimensional, spectrally accurate, integration method by finding the collocation points which minimize the absolute condition number of the integration formula. Although, for the simplicity of the implementation we used the downhill simplex method, which is not an optimal algorithm, we demonstrated that efficient integration method can be derived even for non-continuous weight functions and non-standard domains. More research is needed in analyzing the properties of this integration method and its implementations. In particular,  finding a proper optimization algorithm for obtaining the collocation points. 
 
It should be noted that this integration method is less efficient than the classical Gauss quadrature as an $M+1$ point formula  is accurate only for the $M+1$ polynomials which are used, rather than achieving increased algebraic accuracy, as in the case of Gauss integration. However, it can be used even when no Gauss-type quadrature can be derived.

Numerical examples demonstrated the efficacy of the proposed method.

\newpage

\appendix
%\section{Appendix - L-Shaped Domain - Local Errors - Non-standard Polynomials}

\section{Appendix -  Local Errors for the L-Shaped Domain }

\begin{table}[h]
%\caption{The error in each subdomain, $\Omega_j$ for different degree of non-stabdard polynomial expansions}
\begin{tabular}{| c | c | c | c | c |}
\hline
Element/Degree & 1 & 2 & 3 & 4 \\ \hline
1 & 0.000601801 & 0.0000504182 & 0.00000129901 & 0.0000002883369 \\ \hline
2 & 0.000283446 & 0.0000621132 & 0.00000120905 & 0.000000276578 \\ \hline
3 & 0.000146762 & 0.0000700263 & 0.00000029024 & 0.00000006435543 \\ \hline
4 & 0.000686025 & 0.0000292184 & 0.00000148222 & 0.000000339858 \\ \hline
5 & 0.000545596 & 0.0000458312 & 0.00000117826 & 0.000000251762 \\ \hline
6 & 0.000283395 & 0.0000621132 & 0.000000605555 & 0.000000405396 \\ \hline
7 & 0.00051755 & 0.00000994984 & 0.00000108778 & 0.0000000131404 \\ \hline
8 & 0.000686025 & 0.0000292184 & 0.00000148222 & 0.000000388262 \\ \hline
9 & 0.000601804 & 0.0000504182 & 0.00000129902 & 0.000000221688 \\ \hline
10 & 0.00060181 & 0.0000504181 & 0.0000012872 & 0.000000221688 \\ \hline
11 & 0.000686025 & 0.0000292172 & 0.0000014835 & 0.000000388264 \\ \hline
12 & 0.00051733 & 0.0000099498 & 0.0000010866 & 0.00000001314 \\ \hline
13 & 0.000611872 & 0.0000013364 & 0.00000013318 & 0.0000000080571 \\ \hline
14 & 0.000836335 & 0.0000157323 & 0.00000180754 & 0.000000444328 \\ \hline
15 & 0.000878226 & 0.0000372123 & 0.00000189706 & 0.000000435895 \\ \hline
16 & 0.00028344 & 0.0000621132 & 0.000000601375 & 0.000000405395 \\ \hline
17 & 0.000545587 & 0.0000458312 & 0.00000117428 & 0.000000251762 \\ \hline
18 & 0.000686024 & 0.0000292184 & 0.0000014833 & 0.000000339858 \\ \hline
19 & 0.000836335 & 0.0000157323 & 0.0000018075 & 0.000000324344 \\ \hline
20 & 0.00102394 & 0.00000221364 & 0.00000221263 & 0.000000523526 \\ \hline
\end{tabular}
\end{table}
\newpage
\begin{table}[h]
\begin{tabular}{| c | c | c | c | c |}
\hline
Element/Degree & 1 & 2 & 3 & 4 \\ \hline
21 & 0.00113068 & 0.0000211816 & 0.0000024426 & 0.000000557253 \\ \hline
22 & 0.000146762 & 0.0000700265 & 0.00000029033 & 0.000000064355 \\ \hline
23 & 0.000283446 & 0.0000621111 & 0.0000121903 & 0.000000276578 \\ \hline
24 & 0.000836335 & 0.0000157323 & 0.0000018075 & 0.000000444328 \\ \hline
25 & 0.000878226 & 0.0000372123 & 0.00000189706 & 0.000000453976 \\ \hline
26 & 0.00113068 & 0.0000211816 & 0.0000024426 & 0.000000658725 \\ \hline
27 & 0.00131677 & 0.0000028449 & 0.0000028445 & 0.0000003862272 \\ \hline
\end{tabular}
\vspace{0.5cm}
\caption{The error in each subdomain, $\Omega_j$ for different degree of non-stabdard polynomials expansions.}
\end{table}

\newpage
\newpage

%\section{Appendix - L-Shaped Domain - Local Errors - Legendre-type Polynomials}
\begin{table}[h]
\begin{tabular}{| c | c | c | c | c |}
\hline
Element/Degree & 1 & 2 & 3 & 4 \\ \hline
1 & 0.00136001 & 0.0000498848 & 0.00000950256 & 0.0000021083 \\ \hline
2 & 0.00116796 & 0.0000847799 & 0.00000402961 & 0.000000894066 \\ \hline
3 & 0.000644328 & 0.000103328 & 0.0000017648 & 0.000000391563 \\ \hline
4 & 0.00177759 & 0.000151793 & 0.0000335481 & 0.000007443461 \\ \hline
5 & 0.00168033 & 0.0000926857 & 0.00000833637 & 0.000001849625 \\ \hline
6 & 0.00116796 & 0.0000847799 & 0.00000402961 & 0.000003597997 \\ \hline
7 & 0.00417917 & 0.00119975 & 0.000653253 & 0.000380225 \\ \hline
8 & 0.00177759 & 0.000151793 & 0.0000335481 & 0.0000074434 \\ \hline
9 & 0.00136001 & 0.0000498848 & 0.0000095025 & 0.00000210837 \\ \hline
10 & 0.00136001 & 0.0000498848 & 0.0000095025 & 0.00000210838 \\ \hline
11 & 0.00177759 & 0.000151793 & 0.0000335481 & 00.00000744342 \\ \hline
12 & 0.00417917 & 0.00119975 & 0.000653253 & 0.000380225 \\ \hline
13 & 0.00295143 & 0.00128745 & 0.000676481 & 0.000384855 \\ \hline
14 & 0.000810902 & 0.000195597 & 0.0000354234 & 0.00000847482 \\ \hline
15 & 0.0011291 & 0.0000536411 & 0.0000096454 & 0.00000230760 \\ \hline
16 & 0.00116796 & 0.0000847799 & 0.0000040296 & 0.00000359799 \\ \hline
17 & 0.00168033 & 0.0000926857 & 0.0000083363 & 0.00000184962 \\ \hline
18 & 0.00177759 & 0.000151793 & 0.0000335481 & 0.0000074434 \\ \hline
19 & 0.000810902 & 0.000195597 & 0.0000354234 & 0.0000084748 \\ \hline
20 & 0.000494571 & 0.0000997282 & 0.0000156367 & 0.0000037409 \\ \hline
\end{tabular}
\end{table}
\newpage
\begin{table}[h]
\begin{tabular}{| c | c | c | c | c |}
\hline
Element/Degree & 1 & 2 & 3 & 4 \\ \hline
21 & 0.000924513 & 0.0000530707 & 0.0000053995 & 0.00000123179 \\ \hline
22 & 0.000644328 & 0.000103328 & 0.0000017648 & 0.00000039156 \\ \hline
23 & 0.00116796 & 0.0000847799 & 0.0000040296 & 0.000000894 \\ \hline
24 & 0.00136001 & 0.0000498848 & 0.000009502 & 0.000002108 \\ \hline
25 & 0.0011291 & 0.0000536411 & 0.0000096454 & 0.0000023076 \\ \hline
26 & 0.00092451 & 0.0000530707 & 0.0000053995 & 0.000001231895 \\ \hline
27 & 0.00104132 & 0.0000447484 & 0.0000024735 & 0.00000059165 \\ \hline
\end{tabular}
\vspace{0.5cm}
\caption{The error in each subdomain, $\Omega_j$ for different degree of Legendre-type polynomials expansions.}
\end{table}

\newpage
\bibliographystyle{amsplain}

\bibliography{GPC_ref}

%\bibliographystyle{amsplain}
%\bibliography{xbib}

%% not want to use model1b-num-names.bst.


%  \begin{thebibliography}{00}

\end{document}